 \pdfoutput=1
\newif\ifarxiv
\arxivtrue 

\documentclass[journal]{IEEEtran} 

\usepackage{amsmath}
\usepackage{amsthm}
\usepackage{amssymb}
\usepackage{graphicx} 
\usepackage{caption}
\usepackage{subcaption}
\usepackage{tikz}
\usepackage{diagbox}
\usepackage{xcolor}
\usepackage{xfrac}
\usepackage{booktabs}
\usepackage{longtable}
\usepackage{multirow}
\usepackage{multicol}
\usepackage{pgfplots}
\usetikzlibrary{arrows.meta, calc, decorations.pathreplacing}
\pgfplotsset{compat=1.17}

\def\ee{0.2} 

\def\q{0.87*1.5*pi} 

\usepackage{cite}
\usepackage{lipsum}

\usepackage{color}

\usepackage{ifthen}
\newboolean{com}
\setboolean{com}{false}

\usepackage{algorithm} 
\usepackage{algpseudocode} 

\newtheorem{theorem}{Theorem}

\newtheorem{lemma}{Lemma}
\newtheorem{definition}{Definition}
\newtheorem{proposition}{Proposition}

\newtheorem{ex}{Example}

\newcommand{\mcl}[1]{\mathcal{#1}}

\newcommand{\R}{\mathbb{R}}
\newcommand{\x}{\mathbf{x}}

\newcommand{\N}{\mathbb{N}}

\newcommand{\norm}[1]{\lVert{#1}\rVert}

\newcommand{\eps}{\varepsilon}

\DeclareMathOperator{\sign}{sign}
\DeclareMathOperator*{\esssup}{ess\,sup}

\title{\LARGE \bf
	Bounding the Error of Value Functions in Sobolev Norm Yields Bounds on Suboptimality of Controller Performance
}

\author{Morgan Jones,%
	\thanks{M. Jones is with the School of Electrical and Electronic Engineering,
		The University of Sheffield, e-mail:   {\tt \scriptsize morgan.jones@sheffield.ac.uk} } 
	Matthew M. Peet
	\thanks{M. Peet is with the School for the Engineering of Matter, Transport and Energy, Arizona State University, e-mail: {\tt \small mpeet@asu.edu } }
}
\begin{document}

\maketitle
\thispagestyle{empty}
\pagestyle{empty}

\begin{abstract}
	Optimal feedback controllers for nonlinear systems can be derived by solving the Hamilton-Jacobi-Bellman (HJB) equation. However, because the HJB equation is a nonlinear partial differential equation, numerical methods typically provide only approximate solutions. While numerical error bounds on approximate HJB solutions are often available, these bounds do not necessarily translate into guarantees on the suboptimality of the resulting controllers. In this paper, {we establish that the suboptimality of the resulting controller is bounded by the $L^\infty$ norm of the HJB residual, which is, in turn, bounded by numerical error in the value function as measured in the Sobolev $W^{1,\infty}$ norm. This implies that convergence of value functions in $W^{1,\infty}$ result in controllers that yield a cost that is arbitrarily close to the true minimum. In contrast, we demonstrate that such guarantees do not hold when the value function error is measured in weaker norms, such as the Sobolev $W^{1,p}$ norm for finite $p$.} These results apply to systems governed by Lipschitz continuous dynamics over a finite time horizon with compact input space. 
\end{abstract}
\vspace{-0.5cm}
\section{Introduction}
\vspace{-0.1cm}
{Consider} finite horizon Optimal Control Problems (OCPs) in continuous time. The objective of an OCP is to determine an input that guides the system dynamics to minimize a given cost function over a specified time horizon. Arguably one of the most popular approaches for solving OCPs is Dynamic Programming (DP), which requires solution of the Hamilton-Jacobi-Bellman (HJB) Partial Differential Equation (PDE). The solution of the HJB PDE, known as the Value Function (VF), is then used to design an optimal controller.

Unfortunately, solving the HJB PDE is notoriously difficult. In general, the HJB PDE does not typically possess a classical smooth solution, and it wasn't until the 1980s that a suitable alternative, known as the viscosity solution, was introduced~\cite{crandall1983viscosity}. To approximate this viscosity solution/value function, several classes of numerical methods have been developed. Such methods include the (mixed) finite elements approach~\cite{gallistl2020mixed} for which one may obtain an approximate VF with error bound in terms of the first order mixed $L^2$ norm (assuming the Cordes condition holds). Other examples of this class of methods includes discretization approaches such as~\cite{kalise2018polynomial} or~\cite{achdou2008homogenization}, which yield an approximate VF with an $L^{\infty}$ error bound --  a bound which convergse as the level of discretization increases. {Semi-Lagrangian (SL) methods~\cite{falcone2013semi} approximate the solution by discretizing the OCP directly, rather than using finite difference approximations of the HJB PDE derivatives. SL methods synthesize feedback controls via a discrete-time lookahead procedure, distinct from the standard instantaneous feedback policy used in continuous control considered in this paper, see Eq.~\eqref{eq: controller synth from VF}. While bounds on the local one-step optimality gap can be derived for SL schemes under certain assumptions (see Chapter 8 in~\cite{falcone2013semi}), these cannot be used in the context of continuous-time controllers/policies defined over the entire trajectory (as in Theorem~\ref{thm: performance bounds} from this paper).} Alternative notable numerical HJB PDE solution schemes include the method of characteristics~\cite{liberzon2011calculus}, max-plus methods~\cite{mceneaney2007curse}, grid based methods~\cite{kang2017mitigating, kunisch2004hjb}, Neural Network methods~\cite{abu2005nearly,tassa2007least}, policy iteration~\cite{vrabie2009neural,meng2024physics} and value iteration~\cite{xiao2023convergence}.

An important emerging field for approximately solving the HJB PDE is that of Physics Informed Neural Networks (PINNs)~\cite{liu2022physics,shilova2023revisiting,mukherjee2023bridging}. In this class of methods, a neural network is used to parameterize solutions of HJB PDE and the residual, as evaluated at collocation (training data) points, is minimized. However, as shown in~\cite{wang20222}, minimizing the Euclidean norm of the HJB residual may result in PINN methods failing to converge to the VF. {We show, minimizing the HJB residual in the $L^\infty$ norm does yield controllers with arbitrarily good performance, however, we also show that even if a method converges to the VF in some sense, such as in the $W^{1,p}$ norm for finite $p<\infty$, this may result in a controller with poor performance no matter how tight the VF approximation is, unless $p=\infty$. This is illustrated in Example~\ref{ex:lib} from Section~\ref{sec: numerical ex}.}

While the use of approximate VFs to construct controllers has been studied, such controllers typically: apply only to OCPs with specific structure (e.g. affine in input~\cite{ribeiro2020control}); do not have associated performance bounds; and/or have conservative assumptions like the differentiability of the VF~\cite{jiang2015global,abu2005nearly,baldi2012scalable,baldi2015piecewise,zhu2017policy}. The work of~\cite{leong2014optimal} does provide a performance bound but requires the assumption that there exists a transformation that linearizes the HJB PDE.  Alternatively~\cite{cunis2020sum} proposes a bilinear SOS optimization framework without any convergence guarantees. The SOS and moment based approaches~\cite{korda2016controller,zhao2017control,chen2020optimal} yield controllers that converge to optimal performance under specific assumptions, however, this convergence is not related to VF approximation.
By contrast, the results of this paper provide a general performance bound independent of the approach used to approximate the VF.


%

For discrete time OCPs over infinite time horizons it has been shown in Proposition 6.1 from~\cite{bertsekas1995neuro} (and also presented in~\cite{singh1994upper}) that if a function $J$, is $\eps$-close to the true VF under the uniform norm, then the sub-optimality of the controller synthesized using $J$ is bounded by ${2 \alpha \eps}/{(1-\alpha)}$, where $\alpha>0$ is the discount factor of the OCP. In this paper we seek to derive an analogous result for the optimality gap of controllers synthesized using approximate VFs for continuous time OCPs over finite time horizons. The main result of the paper, given in Theorem~\ref{thm: performance bounds}, shows that controller performance is bounded by $C\eps$ where $\eps>0$ is the VF approximation error in the $W^{1, \infty}$ norm and $C>0$ is a known constant depending on the problem data.

\textbf{Notation:} For $x \in \R^n$ we denote the Euclidean norm by $\|x\|_2=\sqrt{\sum_{i=1}^n x_i^2 }$. For $r>0$ and $x \in \R^n$ we denote the open ball $B_r(x):=\{y \in \R^n: \|x-y\|_2<r\}$ and closed ball $\bar{B}_r(x):=\{y \in \R^n: \|x-y\|_2 \le r\}$. 
For $p \in [1,\infty)$ we denote the set of $p$-integrable functions by $L^p(\Omega):=\{f:\Omega \to \R \text{ measurable }: \int_{\Omega}|f|^p<\infty    \}$, in the case $p = \infty$ we denote $L^\infty(\Omega):=\{f:\Omega \to \R \text{ measurable }: \esssup_{x \in \Omega}|f(x)| < \infty \}$.  We use the notation a.e to mean almost everywhere, meaning the set in which a property does not hold has measure zero. For $\alpha \in \N^n$ we denote the partial derivative $D^\alpha f(x):= \Pi_{i=1}^{n} \frac{\partial^{\alpha_i} f}{\partial x_i^{\alpha_i}} (x)$. Let $C(\Omega,\Theta)$, $AC(\Omega,\Theta)$ and $Lip(\Omega,\Theta)$ be the set of continuous, absolutely continuous and uniformly Lipschitz continuous 
functions respectively with domain
$\Omega \subseteq \R^n$ and image $\Theta \subseteq \R^m$. We denote the set  of $i$ continuously differentiable functions by $C^i(\Omega,\Theta):=\{f \in C(\Omega,\Theta): D^\alpha f \in C(\Omega, \Theta) \text{ } \text{ for all } \alpha \in \N^n,\; \norm{\alpha}_{1} \le i\}$. For $k \in \N$ and $1 \le p \le \infty$ we denote the Sobolev space of functions with weak derivatives by $W^{k,p}(\Omega):=\{u\in L^p(\Omega): D^\alpha u \in L^p(\Omega) \text{ for all } \alpha \in \N^n,\; \norm{\alpha}_\infty \le k \}$. For $u \in W^{k,p}(\Omega)$ we denote the Sobolev norm 
$\|u\|_{W^{k,p}(\Omega)}:= \begin{cases}
	\left( \sum_{|\alpha| \le k} \int_\Omega |D^\alpha u(x)|^p dx \right)^{\frac{1}{p}} \text{ if } 1 \le p < \infty\\
	\sum_{|\alpha| \le k} \esssup_{ x \in \Omega  } \{|D^\alpha u(x) |\} \text{ if } p= \infty.
\end{cases}$ 
In the case $k=0$ we have $W^{0,p}(\Omega)=L^p(\Omega)$ and thus we use the notation $\| \cdot \|_{L^p(\Omega)} :=\| \cdot \|_{W^{0,p}(\Omega)} $. For $V \in C^1(\R \times \R^n , \R)$ we denote $\nabla_x V:= (\frac{\partial V}{\partial x_2},....,\frac{\partial V}{\partial x_{n+1}})^{\top}$. {For a function $u \in C(\Omega, \R)$ and $x \in \Omega \subset \R^N$ we define the super/sub differential sets 
\vspace{-0.4cm}\begin{align*}
	D^+ \hspace{-0.05cm}  u(x) & \hspace{-0.1cm}  := \hspace{-0.1cm}  \left\{ \hspace{-0.05cm}  p \in \mathbb{R}^N\hspace{-0.1cm}  : \hspace{-0.1cm}  \limsup_{y \to x, y \in \Omega} \hspace{-0.1cm}  \frac{u(y) - u(x) - p^\top (y - x)}{|y - x|} \hspace{-0.05cm}  \le  \hspace{-0.05cm}  0 \hspace{-0.05cm}  \right\}. \\
	D^- \hspace{-0.05cm}  u(x) & \hspace{-0.1cm}  := \hspace{-0.1cm}  \left\{ \hspace{-0.05cm}  p \in \mathbb{R}^N \hspace{-0.1cm}  : \hspace{-0.1cm} \liminf_{y \to x, y \in \Omega} \hspace{-0.1cm}  \frac{u(y) - u(x) - p^\top (y - x)}{|y - x|} \hspace{-0.05cm}  \ge \hspace{-0.05cm}  0 \hspace{-0.05cm}  \right\}.
\end{align*}
We denote the union of super/sub differential sets  as $D^{\pm} u (x)=D^+u(x) \cup D^-u(x)$.
For $u \in \text{Lip}(\Omega, \mathbb{R})$, we define the \textit{Clarke generalized gradient} at $x \in \Omega$ as
\vspace{-0.3cm}\begin{align} \nonumber
	\partial_C u(x) := \bigg\{ p \in \mathbb{R}^n :   p^\top  v  \le \limsup_{y \to x, h \downarrow 0} & \frac{u(y+hv) - u(y)}{h}\\ \nonumber 
	&  \text{ for }v \in \mathbb{R}^n \bigg\}.
\end{align}
For functions $J$ over $\mathbb{R} \times \mathbb{R}^n$, elements $(p,q) \in D^\pm J(t,x)$ or $(p,q) \in \partial_C J(t,x)$ explicitly separate the time differential $p \in \mathbb{R}$ from the spatial differential $q \in \mathbb{R}^n$.}
\vspace{-0.4cm}
\section{{Optimal Control Problems}}\label{sec:prob-form} 
Consider the following family of Optimal Control Problems (OCPs), each initialized by $(t_0,x_0) \in [0,T] \times \R^n$,
\vspace{-0.05cm}
\begin{align}  \nonumber 
	&V^*(t_0,x_0):=\inf_{ u}   \int_{t_0}^{T} \hspace{-0.2cm} c(t,x(t),  u(t) ) dt + g(x(T)) \text{ subject to:}  \\ \label{opt: optimal control probelm}
	& x(t_0)=x_0, \quad  \dot{x}(t) = f(t,x(t),  u(t)) \text{ for } t \in [t_0,T],\\ \nonumber
	&  u:[t_0,T] \to \R^m \text{ is a measurable function,} \\ \nonumber
	& u(t) \in U \text{ for all } t \in [t_0,T],
\end{align}
where $c: \R \times \R^n \times \R^m \to \R$, $g: \R^n \to \R$, $f:\R \times \R^n \times \R^m \to \R^n$, $U \subset \R^m$, and $T>0$. Each OCP~\eqref{opt: optimal control probelm} is associated with a tuple $\{c,g,f,U,T\}$. {A measurable function, $u:[t_0,T] \to \R^m$, is an \textit{admissible input} for OCP~\eqref{opt: optimal control probelm} if $u(t) \in U$.} We refer to $V^*$, defined in Eq.~\eqref{opt: optimal control probelm}, as the Value Function (VF) of the OCP. Throughout this paper we consider the following standard class of OCPs for which it can be shown that the VF is Lipschitz continuous and solution maps exist.
\vspace{-0.1cm}
\begin{definition} \label{ass: conditions on OCP}
	We say OCP~\eqref{opt: optimal control probelm} corresponding to tuple $\{c,g,f,U,T\}$ is class $L$ if the following hold:
	\begin{enumerate}
		\item $U \subset \R^m$ is a compact set.
		\item $f$, $c$ and $g$ are {uniformly Lipschitz continuous in the state variable and uniformly continuous in the time and input variable.} Specifically, there exist $\alpha_f, \alpha_c,\alpha_g>0$ such that  $\|f(t,x,u)-f(t,y,u)\|_2 \le \alpha_f \|x-y\|_2$, $\|c(t,x,u)-c(t,y,u)\|_2 \le \alpha_c \|x-y\|_2$ and $\|g(x)-g(y)\|_2 \le \alpha_g \|x-y\|_2$  for all $x,y \in \R^n$, $t \in [0,T]$ and $u \in U$.
		\item $f$, $c$ and $g$ are bounded by $\mcl{O}(\|x\|_2)$. Specifically, there exists $\beta_f,\beta_c,\beta_g>0$ such that $\|f(t,x,u)\|_2 \le \beta_f (1+\|x\|_2) $, $\|c(t,x,u)\|_2 \le \beta_c (1+\|x\|_2)$ and $\|g(x)\|_2 \le \beta_g (1+\|x\|_2)$ for all $x \in \R^n$, $t \in [0,T]$ and $u \in U$.
	\end{enumerate}
\end{definition}

{\noindent \textbf{Solution Map.} Under the first two items of Definition~\ref{ass: conditions on OCP}, standard global existence and uniqueness theorems (see, e.g., \cite[Thm.~3.2]{Khalil2002}) guarantee that for any initial condition $(t_0,x_0)$ and any admissible input, the initial value problem $\dot{x}(t) = f(t,x(t), u(t)), \quad x(t_0)=x_0$
admits a unique absolutely continuous solution defined on the interval $[t_0, T]$. We denote this unique trajectory by $x_{u,x_0}(\cdot)$, and refer to the mapping $(u, x_0) \mapsto x_{u,x_0}$ as the {solution map} of the system.}

\begin{definition}[Admissible \& Optimal Controller] \label{def: feedback controller}
	{ Given an OCP~\eqref{opt: optimal control probelm}, we say a measurable function $\pi(t,x)$ is an \textbf{admissible controller} for the OCP if, for any $({t_0},x_0) \in [0,T]\times \R^n$, 
		\vspace{-0.2cm} \[ 
		\dot{x}(t) = f(t,x(t), \pi(t,x(t))), \quad x(t_0)=x_0, 
		\]
		has a unique solution and $\pi(t,x(t)) \in U$ for all $t \in [t_0,T]$. {For $t_0=0$}, we denote the corresponding admissible input  as $u_{\pi,x_0}(t)=\pi(t,x(t))$. Furthermore, an admissible controller, $\pi$, is \textbf{optimal} if $u_{\pi,x_0}$ solves OCP~\eqref{opt: optimal control probelm} {with $t_0=0$}.}
\end{definition}

{The Dynamic Programming (DP) method allows one to obtain an optimal controller by solving HJB PDE:}
\vspace{-0.2cm}
\begin{align}  \nonumber
	&\frac{\partial}{\partial t} V(t,x) + \inf_{u \in U} \left\{ {H(t,x,u, \nabla_x V(t,x))} \right\} = 0 \\
	& \hspace{3.75cm} \text{ for all } (t,x) \in [0,T] \times \R^n,   \label{eqn: general HJB PDE}\\ \nonumber
	& V(T,x)= g(x) \quad \text{ for all } x \in \R^n,
\end{align}
{where we define the Hamiltonian as
 \vspace{-0.1cm} \begin{align} \label{eq: hamiltonian}
	H(t,x,u,q):=c(t,x,u) + q^{\top} f(t,x,u).
\end{align}}
{Given a differentiable solution to the HJB PDE~\eqref{eqn: general HJB PDE}, $V^*$, an optimal controller may be constructed by minimizing the Hamiltonian~\cite{liberzon2011calculus}:
\vspace{-0.1cm}	\begin{align} \label{eq: controller synth from VF}
	\pi^*(t,x) \in \arg \inf_{u \in U} H(t,x,u,\nabla_xV^*(t,x)).
\end{align} }

Note, the HJB PDE may not possess a classical differentiable solution. However, it can be shown under mild assumptions that the HJB PDE will always admit a ``viscosity solution'' -- See, e.g.~\cite{crandall1983viscosity}.
Furthermore, for class L (Defn.~\ref{ass: conditions on OCP}) OCPs, it can be shown that the VF (as defined in Eq.~\eqref{opt: optimal control probelm}) is the unique viscosity solution of the HJB PDE~\eqref{eqn: general HJB PDE} and is Lipschitz continuous.

\begin{lemma}[Page 555 \cite{evans2010partial}] \label{lem: unique solution to HJB PDE}
	Suppose tuple $\{c,g,f,U,T\}$ defines an OCP of class L (Defn.~\ref{ass: conditions on OCP}). Then $V^*:  [0,T]  \times \R^n \to \R$, as defined in Eq.~\eqref{opt: optimal control probelm}, is the unique viscosity solution of the HJB PDE~\eqref{eqn: general HJB PDE} and {$V^* \in Lip([0,T] \times \R^n , \R)$.} 
\end{lemma}



{Eq.~\eqref{eq: controller synth from VF} can be used to determine the value of the optimal controller, $\pi^*(t,x)$, at points where $\nabla_xV^*(t,x)$ exists. Furthermore, Lemma~\ref{lem: unique solution to HJB PDE} implies that for class L (Defn.~\ref{ass: conditions on OCP}) OCPs, $V^*$ is Lipschitz continuous and hence differentiable almost everywhere. Unfortunately, however, segments of the optimal closed-loop trajectory may coincide with the measure zero regions of state space in which $V^*$ is not differentiable \ifarxiv (see Example~\ref{ex: non-diff traj} and Fig.~\ref{fig: moving along the ridge} for an illustration) \fi. In such cases, we consider the approach used in~\cite{zhou1993verification} to define sub/super differentials ($D^{\pm}$) of $V^*$ at these points.}


{\begin{definition} \label{def: admis VF}
Suppose OCP~\eqref{opt: optimal control probelm} is class L (Defn.~\ref{ass: conditions on OCP}). We say a function $J \in Lip( [0,T] \times \R^n , \R)$ is an \textbf{admissible candidate value function} of OCP~\eqref{opt: optimal control probelm} if 
\begin{enumerate}
	\item ${D^\pm J(t,x)} = \partial_C J(t,x)$ for all $(t,x) \in [0,T] \times \R^n$.
	\item There exists an admissible controller $\pi_J: [0,T] \times \R^n \to U$ and measurable functions $(p_J(t,x),q_J(t,x)) \in {D^\pm J(t,x)}$ such that 
	\vspace{-0.1cm} \begin{align} \nonumber
		(p_J(t,x),q_J(t,x), &	\pi_J(t,x))\\ \nonumber
		& \in   \arg \hspace{-1cm} \inf_{(p,q,u) \in D^\pm J(t,x) \times U} \hspace{-0.8cm} \{p+ H(t,x,u,q) \}\\ \label{eq: controller synth}
		& \text{ for all } (t,x) \in [0,T] \times \R^n,
	\end{align}
	where $H$ is given in Eq.~\eqref{eq: hamiltonian}.
\end{enumerate}
\end{definition}
The use of sub/super differentials in defining the controller $\pi_J$, as in Definition~\ref{def: admis VF}, will allow us to extend the applicability of the main result given in Theorem~\ref{thm: performance bounds} to OCPs where the VF is only Lipschitz continuous and not continuously differentiable. However, this definition requires that the  union of sub- and superdifferentials of the candidate VF coincide with the Clarke generalised gradient. This condition is automatically satisfied when $J$ is continuously differentiable, since $D^\pm J(t,x)=\partial_C J(t,x)=\{\nabla J(t,x)\}$ (Lemma~1.8 in~\cite{bardi1997optimal}, Theorem~10.8 in~\cite{clarke2013functional}). Examples of non-differentiable admissible candidate value functions include ridges or valleys -- e.g.  $J(x) = |x|$ or $J(x)=-|x|$, for which $D^\pm J(t,x)=\partial_C J(t,x)=[-1,1]$. 

More generally, semiconcave and semiconvex functions are admissible (Prop.~4.7 in~\cite{bardi1997optimal}), a class that includes VFs arising from a wide range of optimal control problems (see Chapter~4.2 in~\cite{bardi1997optimal}). Examples of non-admissible candidate value functions (where $D^\pm J \ne \partial_C J$) include functions with saddle points (e.g., $J(t,x)=|x_1-t|-|x_2-t|$, where $D^\pm J(t,[t,t]^\top)=\emptyset$) and functions with singular oscillations such as $J(x)=x^{3/2}\sin(1/\sqrt{x})$ -- for which $D^\pm J(0)=\{0\}$ but $\partial_C J(0)=[-1,1]$ (see Exercise~13.10 in~\cite{clarke2013functional}).}

\vspace{-0.2cm}
\section{A Performance Bound On Controllers Constructed Using Approximate Value Functions} \label{sec:optimal controller approx}
\ifarxiv
\begin{figure}
\centering
\begin{tikzpicture}[thick, scale=0.8]
	\begin{axis}[
		axis lines = center,
		xlabel = {},
		ylabel = {},
		xtick=\empty,
		ytick=\empty,
		domain=-1*pi:2*pi, 
		samples=100,
		grid=both,
		trig format = rad,
		grid style={dashed, gray!30},
		thick,
		legend style={at={(0,1)}, anchor=north,legend columns=1}, 
		]
		\addplot[
		smooth,
		blue,
		thick,
		]
		{(1 +\ee * (x - 2*pi) / (3*pi)) * sin(x)};
		\addlegendentry{$H(t,x,u, \nabla \textcolor{blue}{V^*}(t,x))$}

\addplot[
smooth,
red,
thick,
]
{(1 - \ee* (x + pi) / (3*pi)) * sin(x)};
\addlegendentry{$H(t,x,u,\nabla \textcolor{red}{J}(t,x))$}

\draw[-, thick, dashed,blue] 
(axis cs:1.5*pi,{(1 +\ee * (1.5*pi - 2*pi) / (3*pi)) * sin(1.5*pi)}) -- 
(axis cs:1.5*pi,0) 
node[above] {\Large \(\textcolor{blue}{u^*}\)};

\draw[-, thick, dashed,red] 
(axis cs:-0.5*pi,{(1 - \ee* (-0.5*pi + pi) / (3*pi)) * sin(-0.5*pi)}) -- 
(axis cs:-0.5*pi,0) 
node[above] {\Large \(\textcolor{red}{{u_J}}\)};

\draw[<->, thick] 
(axis cs:\q,{(1 +\ee * (\q - 2*pi) / (3*pi)) * sin(\q)}) -- 
(axis cs:\q,{(1 - \ee* (\q + pi) / (3*pi)) * sin(\q)}) 
node[right,midway] {\Large \(\epsilon\)};

\end{axis}
\end{tikzpicture}
\vspace{-5pt}
\caption{\footnotesize
Illustration of control policy sensitivity to VF approximation error. Although the Hamiltonian curves derived from the optimal VF $V^*$ (blue) and the approximate VF $J$ (red) are close ($\epsilon$-bounded), their respective global minima, $u^*$ and $u_J$, are significantly different. This demonstrates that $V^* \approx J$ does not necessarily imply $u_J \approx u^*$ for a fixed time and state.} \label{fig: close VFs}
\end{figure}
\fi

%

{Having defined classes of admissible candidate VFs, we may now proceed to the main result of the paper. Of course, given an OCP, if a viscosity solution, $V^*$ to the HJB PDE~\eqref{eqn: general HJB PDE} is known, an admissible controller can be verified as optimal using \cite[Thm.~3.5]{zhou1993verification}}. However, finding a HJB PDE {viscosity solution} is challenging and hence we consider the case where a function, $J$, is obtained which is close to $V^*$ with respect to some norm, $X$ -- i.e. $\norm{J-V^*}_X\le \epsilon$, but does not necessarily satisfy the HJB PDE. This naturally raises the question: how well does a controller constructed from an approximate VF perform and does the choice of approximation norm matter?

\ifarxiv
As established in this section, the suboptimality of a controller obtained from an approximate VF depends strongly on the norm in which the candidate VF approximates the true VF. For illustration, consider Figure~\ref{fig: close VFs}, which shows that two candidate VFs ($V^*$ and $J$) may be close in the $L_\infty$ norm, {yet the points at which their Hamiltonians~\eqref{eq: hamiltonian} attain their minima ($u^*(t,x)$ and $u_J(t,x)$) may be arbitrarily far apart.}  \fi

To formalize this analysis, then, we first define the optimality gap (loss) of an initial condition, $x_0$ and admissible input, $u$, associated with an OCP~\eqref{opt: optimal control probelm} {with $t_0=0$}: 
\vspace{-0.1cm} \begin{align}\label{eqn:L}
& L( u,x_0)  :=\\ \nonumber
& \int_{0}^{T} c(s,x_{ u, x_0}(s),   u (s)) ds  + g(x_{u, x_0}(T))   -V^*(0,x_0),
\end{align} 
where $V^*$ is the VF of the OCP. Clearly, $L(u,x_0 ) \ge 0$ for any admissible $u$.

{Given any {admissible candidate value function}, $J $, we may define an associated admissible controller, $\pi_J$, as in Eq.~\eqref{eq: controller synth} (or equivalently Eq.~\eqref{eq: controller synth from VF} if $J$ is differentiable). Our goal is then to bound the optimality gap $L(u_{\pi_J,x_0},x_0 )$ in terms of $\|V^*-J\|_{W^{1,\infty}([0,T]\times \R^n)}$, where $V^*$ is the true VF of the OCP. This result can be stated as follows (with proof deferred).}

\begin{theorem} \label{thm: performance bounds}
Suppose the OCP~\eqref{opt: optimal control probelm} corresponding to tuple $\{c,g,f,U,T\}$ is of class $L$ (Defn.~\ref{ass: conditions on OCP}), and denote the associated value function by $V^*$. {Let $J$ be an {admissible candidate value function} (Definition~\ref{def: admis VF})}, with associated admissible input $u_{\pi_J,x_0}$. Then for any $x_0 \in \R^n$, we have
\vspace{-0.0cm}
\begin{align} \label{ineq:loss functin}
& L( u_{\pi_J,x_0},x_0) \le  C \|J- V^*\|_{W^{1,\infty}([0,T] \times B_R(0))}, \\ \label{eqn: C}
& 	\vspace{-0.2cm} \text{where } C:= 2\max \left\{1,T, T {\beta_f(1+R)} \right\},
\end{align} 
$L$ (the loss) is defined in Eq.~\eqref{eqn:L}, $R:=(1+\|x_0\|_2)e^{\beta_f T}-1$, and $\beta_f$ is as in Defn.~\ref{ass: conditions on OCP}.
\end{theorem}

{To establish Theorem~\ref{thm: performance bounds}, we first show that for any candidate value function there is an associated OCP for which the resulting admissible controller is optimal.}

{\begin{lemma} \label{lem: J solves modified OCP}
Consider OCP~\eqref{opt: optimal control probelm} with tuple $\{c,{g},f,U,T\}$. For any {admissible candidate value function} $J $ {(Definition~\ref{def: admis VF})}, suppose $\pi_J$ satisfies Eq.~\eqref{eq: controller synth}. Then $\pi_J$ is optimal for the OCP defined by the tuple $\{\tilde{c},\tilde{g},f,U,T\}$, where
\begin{align} \label{eq: modified OCP def}
& \tilde{c}(t,x,u)  :=c(t,x,u)-E_J(t,x), \quad   \tilde{g}(x)  :=J(T,x),\\ \nonumber
&	E_J(t,x):= \inf_{(p,q,u)\in D^\pm J(t,x) \times U} \{p + c(t,x,u) + q^{\top} f(t,x,u)\}.
\end{align}
\end{lemma}
\begin{proof}

For any admissible input $u:[0,T] \to U$, let us denote the associated cost as
$\tilde{W}_u(t_0,x_0):= \int_{t_0}^{T} \tilde{c}(t,x_{u,x_0}(t),  u(t) ) dt + \tilde{g}(x_{u,x_0}(T)) $. We show $\pi_J$ is optimal by showing that
\begin{enumerate}
\item $J(0,x_0)\le \tilde{W}_u(0,x_0)$ for any {admissible input}, $u$.
\item $J(0,x_0)=\tilde{W}_{u_{\pi_J,x_0}}(0,x_0)$.
\end{enumerate} 

To begin, consider $(\pi_J, p_J,q_J)$ satisfying Eq.~\eqref{eq: controller synth}. Then $(p_J(t,x),q_J(t,x)) \in D^\pm J(t,x)$ for all $(t,x) \in [0,T] \times \R^n$ and we have
\begin{align} \label{pfeq: HJB PDE}
&	p_J(t,x)+\inf_{u \in U} \{\tilde{c}(t,x,u)+q_J(t,x)^\top f(t,x,u) \}\\ \nonumber 
&	=\inf_{u \in U} \hspace{-0.05cm} \{p_J(t,x)+c(t,x,u)+q_J(t,x)^\top f(t,x,u) \} \hspace{-0.05cm}  - \hspace{-0.05cm} E_J(t,x)\\ \nonumber 
&	=\inf_{(p,q,u)\in D^\pm J(t,x) \times U} \hspace{-0.5cm} \{ p + c(t,x,u) + q^\top f(t,x,u)  \}-E_J(t,x)\\
&	=E_J(t,x)-E_J(t,x)=0.\notag
\end{align}

%


Now, for any admissible input, $u$,  $x_{u,x_0}$ is absolutely continuous which implies by Lemma~\ref{lem:chain_rule} ({found in the appendix}) that  $\frac{d}{dt} J(t,x_{u,x_0}(t))= p + q^\top f(t,x_{u,x_0}(t),u(t))$ a.e for $t \in [0,T]$ and for any $(p,q) \in D^\pm J(t,x_{u,x_0}(t))$. Since $(p_J(t,x_{u,x_0}(t)),q_J(t,x_{u,x_0}(t))) \in D^\pm J(t,x_{u,x_0}(t))$, we have
\begin{align} \label{pfeq:chain rule}
&	\frac{d}{dt} J(t,x_{u,x_0}(t)) \\ \nonumber 
& = p_J(t,x_{u,x_0}(t)) + q_J(t,x_{u,x_0}(t))^\top f(t,x_{u,x_0}(t),u(t))\\ \nonumber
& = p_J(t,x_{u,x_0}(t)) + \tilde{c}(t,x_{u,x_0}(t),u(t)) \\\notag
&+q_J(t,x_{u,x_0}(t))^\top f(t,x_{u,x_0}(t),u(t))-\tilde{c}(t,x_{u,x_0}(t),u(t))\\ \nonumber
& \ge  p_J(t,x_{u,x_0}(t))-\tilde{c}(t,x_{u,x_0}(t),u(t))\\\notag
&+\inf_{v \in U} \{\tilde{c}(t,x_{u,x_0}(t),u)+q_J(t,x_{u,x_0}(t))^\top f(t,x_{u,x_0}(t),v) \}\\\notag
&= - \tilde{c}(t,x_{u,x_0}(t),u(t)) \text{ a.e } t \in [0,T].\notag
\end{align}
Where the last equality follows from Eq.~\eqref{pfeq: HJB PDE}.
%
%
Integrating Eq.~\eqref{pfeq:chain rule} over $[0,T]$ and using $J(T,x)=\tilde{g}(x)$ we deduce $J(0,x_0)\le \tilde{W}_u(0,x_0)$ for any admissible input.

Now consider the admissible input, $u_{\pi_J,x_0}$, obtained from the candidate value function $J$ and associated controller, $\pi_J$. To reduce notational complexity we denote the trajectory associated with this input, $x_{u_{\pi_J,x_0},x_0}$, as $x_J$.  Then, by Eq.~\eqref{pfeq: HJB PDE} and since it follows from Eq.~\eqref{eq: controller synth} that the infimum in Eq.~\eqref{pfeq: HJB PDE} is attained by $u_{\pi_J,x_0}$ we have
\begin{align*}
&p_J(t,x_J(t)) + \tilde{c}(t,x_J(t),u_{\pi_J,x_0}(t)) \\ 
\nonumber & \qquad +q_J(t,x_J(t))^\top f(t,x_J(t),u_{\pi_J,x_0}(t))\\ \nonumber
& =  p_J(t,x_J(t))+\\ \nonumber
&\quad \inf_{v \in U} \{\tilde{c}(t,x_J(t),v)+q_J(t,x_J(t))^\top f(t,x_J(t),v) \}=0,
\end{align*}
and hence as per Eq.~\eqref{pfeq:chain rule}, $\frac{d}{dt} J(t,x_J(t))=- \tilde{c}(t,x_J(t),u_{\pi_J,x_0}(t))$. We conclude by integrating that $J(0,x_0)= \tilde{W}_{u_{\pi_J,x_0}}(0,x_0)$ -- establishing that the associated control $\pi_J$ is optimal.
\end{proof}
{Before deriving our controller performance bound (given in Theorem~\ref{thm: performance bounds}) we next prove $E_J$, as defined in Eq.~\eqref{eq: modified OCP def}, is lower-semi continuous for any admissible candidate value function, $J$ -- allowing us to apply Lemma~\ref{lem: ess sup bound} {(found in the appendix)}.}
\begin{lemma} \label{lem: E LSC}
Consider OCP~\eqref{opt: optimal control probelm} with tuple $\{c,{g},f,U,T\}$. Let $J$ be an {admissible candidate value function} {(Definition~\ref{def: admis VF})}, then
$E_J(t,x):= \inf_{(p,q,u)\in D^\pm J(t,x) \times U} \{p + c(t,x,u) + q^{\top} f(t,x,u)\}$
is lower semi-continuous.
\end{lemma}
\begin{proof}
	To show $E_J$ is lower semi-continuous, we show that for any bounded sequence $(t_n,x_n) \to (t,x)$, we have $\displaystyle E_J(t,x) \le \liminf_{n \to \infty} E_J(t_n,x_n)$.
	
	First, by definition of $\liminf$ there exists a subsequence $\{ t_{n_k}, x_{n_k} \}_k$ such that $\lim_{k \rightarrow \infty} E_J(t_{n_k}, x_{n_k})=\liminf_{n \to \infty} E_J(t_n,x_n)$. Now, for each $k$ in the subsequence, let $\{p_k, q_k, u_k\}$ be the sequence of associated minimizers so that:
	\vspace{-0.2cm}$$
	E_J(t_{n_k}, x_{n_k}) = p_k + H(t_{n_k}, x_{n_k}, u_k, q_k),
	$$ where $H$ is given in Eq.~\eqref{eq: hamiltonian}.
    Let $L_J$ be a uniform Lipschitz constant for $J\in Lip( [0,T] \times \R^n , \R)$. It follows that $\partial_C J(t,x) \subset \bar{B}_{L_J}(0)$ for $(t,x) \in [0,T] \times \R^n$ (See, e.g.~\cite{clarke2013functional} (Prop 10.5)). Since $(p_k, q_k, u_k) \in D^\pm J(t_{n_k},x_{n_k}) \times U = \partial_C J(t_{n_k},x_{n_k}) \times U \subset \bar{B}_{L_J}(0) \times U$ for all $k$, and both $\bar{B}_{L_J}(0)$ and $U$ are compact, the sequence $\{p_k, q_k, u_k\}$ is uniformly bounded.
We now, using the Bolzano–Weierstrass theorem, extract a convergent subsequence $\{p_{k_j}, q_{k_j}, u_{k_j}\}$ of $\{p_k, q_k, u_k\}$  where $\lim_{j\rightarrow \infty }(p_{k_j}, q_{k_j}, u_{k_j}) \to (p^*, q^*, u^*)$. 
Because $U$ is closed we get $u^* \in U$ and Prop. 10.10 in~\cite{clarke2013functional} implies $(p^*, q^*) \in \partial_C J(t,x)$.

	
	By the continuity of the Hamiltonian $H$:
	\vspace{-0.2cm} \begin{align*}
	\lim_{j \to \infty} (p_{k_j} + H(t_{n_{k_j}}, x_{n_{k_j}}, u_{k_j}, q_{k_j})) = p^* + H(t, x, u^*, q^*).
	\end{align*}
Hence,
\begin{align*}
\liminf_{n \to \infty} E_J(t_n, x_n)&=\lim_{k \rightarrow \infty} E_J(t_{n_k}, x_{n_k})=\lim_{j \rightarrow \infty} E_J(t_{n_{k_j}}, x_{n_{k_j}})\\
&=\lim_{j \to \infty} (p_{k_j} + H(t_{n_{k_j}}, x_{n_{k_j}}, u_{k_j}, q_{k_j}))\\
&=p^* + H(t, x, u^*, q^*).
\end{align*}

	Since $(p^*, q^*, u^*)\in \partial_C J(t,x) \times U = D^\pm J(t,x) \times U$, for any $(t,x)$, by definition of $E_J$,
\[
	E_J(t,x) \le p^* + H(t, x, u^*, q^*) = \liminf_{n \to \infty} E_J(t_n, x_n).
\]
\end{proof}

%
 We now bound the optimality gap of using a controller/policy synthesized from an approximated VF based on the HJB PDE residual error.
\vspace{-0.1cm}\begin{proposition} \label{prop: Loss bounded by HJB res} 
Suppose the OCP~\eqref{opt: optimal control probelm} corresponding to tuple $\{c,g,f,U,T\}$ is of class $L$ (Defn.~\ref{ass: conditions on OCP}). Then for any any {admissible candidate value function} $J $ {(Definition~\ref{def: admis VF})} with associated admissible controller, $\pi_J$ Eq.~\eqref{eq: controller synth}, we have that
\vspace{-0.1cm}
\begin{align} \label{ineq:loss functin res}
&	L( u_{\pi_J,x_0},x_0) \le   2 T\delta_1+2\delta_2, \\
&\delta_1 =\hspace{-3mm} \esssup_{x \in B_R(0),t \in [0,T]}\left|\frac{\partial}{\partial t}J(t,x)  + \inf_{u \in U} \{H(t,x,u, \nabla_x J(t,x)) \}\right|, \notag \\ \nonumber 
& \delta_2= \esssup_{x \in B_R(0)} \left|J(T,x)-g(x) \right|,
\end{align} 
where $L$ (the loss) is defined in Eq.~\eqref{eqn:L}, $R:=(1+\|x_0\|_2)e^{\beta_f T}-1$ and $\beta_f$ is as in Defn.~\ref{ass: conditions on OCP}.
\end{proposition}
\textbf{Note:} $J$ is Lipschitz so the derivatives, $\frac{\partial}{\partial t}J(t,x)$ and $\nabla_x J$, in $\delta_1$ are defined a.e and hence the $\esssup$ is well-defined. }
\begin{proof}
By Lemma~\ref{lem: J solves modified OCP} we have that $u_{\pi_J,x_0}$ solves the OCP associated with tuple $\{\tilde{c},\tilde{g},f,U,T\}$, where $\tilde{c}$ and $\tilde{g}$ are defined in terms of $E_J$ as in Eq.~\eqref{eq: modified OCP def}. Denoting its associated trajectory simply as $x_J$ (rather than $x_{u_{\pi_J,x_0},x_0}$), we have for any admissible input $u: [0,T] \to U$ that
\vspace{-0.15cm}\begin{align} \label{u_J ineq}
& \int_{0}^{T} \tilde{c}( s, x_J(s),  u_{\pi_J,x_0}(s)) ds  + \tilde{g}(x_J(T)) \\ \nonumber
& \le \int_{0}^{T} \tilde{c}(s, x_{ u, x_0}(s),   u(s)) ds  + \tilde{g}(x_{ u, x_0}(T)).
\end{align}
By substituting $\tilde{c}(t,x,u)=c(t,x,u)-E_J(t,x)$ and $\tilde{g}(x)=J(T,x)$  into Inequality~\eqref{u_J ineq} and re-arranging terms, we have
\vspace{-0.1cm}\begin{align} \label{difference in cost}
& \int_{0}^{T} c(  s,x_J(s),  u_{\pi_J,x_0}(s)) ds  + {g}(x_J(T)) \\ \nonumber
&  \qquad - \int_{0}^{T} {c}(s,x_{ u, x_0}(s),   u (s)) ds  - {g}(x_{ u, x_0}(T)) \\ \nonumber
& \le \int_0^{T} E_J(s,x_J(s))  - E_J(s,x_{ u, x_0}(s))  ds \\\nonumber
& \qquad + g(x_J(T)) - J(T,x_J(T))\\ \nonumber
& \qquad + J(T,x_{ u, x_0}(T)) - g(x_{ u, x_0}(T)) \\ \nonumber
& \le T\esssup_{s \in[0,T]}\{ | E_J(s,x_J(s))    - E_J(s,x_{ u, x_0}(s))  | \}\\ \nonumber
& \qquad +  2\sup_{y \in B_R(0)}\{ |g(y) - J(T,y)| \}\\ \nonumber
& \le {2T \hspace{-0.5cm}  \esssup_{(s,y) \in [0,T] \times B_R(0)} \hspace{-0.5cm} \{ |E_J(s,y)| \} +  2 \hspace{-0.05cm}  \esssup_{y \in B_R(0)}\{ |g(y) - J(T,y)| \}.}
\end{align} Here the first inequality is simply Eq.~\eqref{u_J ineq}. {To establish the second inequality, we bound the integral using the monotonicity of the Lebesgue integral (H\"{o}lder's inequality). This operation is well-defined because the mapping $t \mapsto E_J(t,x(t))$ is measurable, which follows automatically from the fact that $E_J$ is lower semi-continuous (Lemma~\ref{lem: E LSC}) and the trajectories are absolutely continuous. Furthermore, the essential supremum is finite because $E_J$ is bounded, a direct consequence of the uniform boundedness of the Clarke generalized derivatives due to the Lipschitz continuity of $J$ (Prop 10.5 \cite{clarke2013functional}). The third inequality relaxes the bound by taking the supremum over a fixed spatial domain rather than the specific trajectories. This is done through the use of Lemma~\ref{lem: ess sup bound} (found in the appendix) and because $x_{u,x_0}(t) $ is contained within the closure of ${B_R}(0)$ for all $t \in [0,T]$ and any admissible input $u:[0,T] \to U$ (since under Defn.~\ref{ass: conditions on OCP} $\|x_{ u, x_0}(t)\|_2 \le (1+\|x_0\|_2)e^{\beta_f T}-1 \text{ for all } t \in [0,T] $ by the Gronwall Bellman inequality, see~\cite{caillau2023algorithmic}). Furthermore, since $V^*$ and $J$ are continuous, their essential supremum over $B_R(0)$ coincides with the standard supremum.}

{Since $J$ is Lipschitz continuous, it is differentiable almost everywhere by Rademacher's theorem, and hence,
\vspace{-0.15cm} \begin{align} \label{pfeq: simplfication of E ae}
& E_J(t,x)=\frac{\partial}{\partial t} J(t,x)  + \inf_{u \in U} H(t,x,u,\nabla_x J(t,x))\\ \nonumber & \hspace{3.5cm} \text{ for } (t,x) \in [0,T] \times B_R(0) \text{ a.e}, \\ \label{pfeq:delta 1 eq}
& \text{ implying } \quad  \delta_1=\esssup_{(s,y) \in [0,T] \times B_R(0)} \hspace{-0.5cm} \{ |E_J(s,y)| \}.
\end{align}}
\ifarxiv Finally, since Inequality~\eqref{difference in cost} holds for all admissible inputs, we have
\vspace{-0.1cm}\begin{align} 
& \sup_{u  \, \text{admissible} } \bigg\{ \int_{0}^{T} c(  s,x_{\pi_J,x_0}(s),  u_{\pi_J,x_0}(s)) ds  + {g}(x_{\pi_J,x_0}(T)) \nonumber \\
& \qquad \qquad  - \int_{0}^{T} {c}(s,x_{ u, x_0}(s),   u (s)) ds  - {g}(x_{ u, x_0}(T)) \bigg\} \notag \\ 
& = \int_{0}^{T} c(  s,x_{\pi_J,x_0}(s),  u_{\pi_J,x_0}(s)) ds  + {g}(x_{\pi_J,x_0}(T)) - V^*(0,x_0)\notag  \\ 
& = L( u_{\pi_J,x_0},x_0)\le 2T \delta_1 +  2   \delta_2.\notag
\end{align}
\else
Finally, since Inequality~\eqref{difference in cost} holds for all admissible $u$, {taking the supremum over $u$ and recognizing that the infimum with respect to admissible inputs of the cost yields the optimal value function $V^*(0,x_0)$, we obtain Eq.~\eqref{ineq:loss functin res} as required.}
\fi
\end{proof}
\vspace{-0.2cm}

  {  We next obtain a bound on the optimality gap by connecting the HJB residual to the Sobolev approximation error between the candidate and true VF. This constitutes the paper's main result, previously stated in Theorem \ref{thm: performance bounds}. }

\vspace{-0.2cm}
\begin{proof}[\textbf{Proof of Theorem~\ref{thm: performance bounds}}]

Throughout this proof we denote the set of points where {a Lipschitz function, $W \in Lip([0,T] \times \R^n, \R)$,} is  differentiable by,
\vspace{-0.2cm} \begin{align*}
&S_{W} \hspace{-0.05cm} := \hspace{-0.05cm} \{(t,x) \hspace{-0.05cm}  \in \hspace{-0.05cm}   [0,T] \hspace{-0.05cm}  \times  \hspace{-0.05cm}  B_R(0) \hspace{-0.05cm}  : \hspace{-0.05cm}   W \hspace{-0.05cm}  \text{ is differentiable at } \hspace{-0.05cm}  (t,x) \}.
\end{align*} \vspace{-0.55cm}
\text{ }\\
From Prop.~\ref{prop: Loss bounded by HJB res}, Eq.~\eqref{pfeq:delta 1 eq} and $V^*(T,x)=g(x)$ it follows that  
\vspace{-0.2cm}	\begin{align} \label{pfeq: loss bound}
{L( u_{\pi_J,x_0},x_0) \le 2T }\hspace{-0.5cm}  & \esssup_{(s,y) \in [0,T] \times B_R(0)} \hspace{-0.5cm} \{ |E_J(s,y)| \} \\ \nonumber 
& +  2 \hspace{-0.05cm}  \esssup_{y \in B_R(0)}\{ |V^*(T,y)- J(T,y)| \}.
\end{align}

We now split the remainder of the proof into three parts. In Part~1, we derive an upper bound for $E_J(s,y) $. In Part~2, we find a lower bound for $ E_J(s,y)$. In Part~3 we combine these to bound $\esssup_{(s,y) \in [0,T] \times B_R(0)} \{ |E_J(s,y)| \}$ in Eq.~\eqref{pfeq: loss bound}
and hence obtain Eq.~\eqref{ineq:loss functin} to complete the proof.

%

\underline{\textbf{Part 1 of Proof:}}
For each $(s,y) \in S_{V^*} {\cap S_{J}}$ consider any
\vspace{-0.2cm} \begin{align*}
k^*_{s,y} \in \arg \inf_{u \in U} \left\{H(s,y,u,\nabla_x V^*(s,y))\right\}. 
\end{align*}
Such a $k^*_{s,y}$ exists for each fixed $(s,y) \in S_{V^*} {\cap S_{J}}$ by the extreme value theorem\footnote{This follows because $U \subset \mathbb{R}^m$ is compact, $c$ and $f$ are continuous, and $\nabla_x V^*$ is independent of $u$ and bounded by Rademacher's Theorem.}.
It follows from Eq.~\eqref{pfeq: simplfication of E ae} that
\vspace{-0.2cm} \begin{align} \label{11}
E_J(s,y) & = \frac{\partial}{\partial s} J(s,y)+ \inf_{u \in U} H(s,y,u, \nabla_x J(s,y))\\ \nonumber
& \le \frac{\partial}{\partial s} J(s,y)+  H(s,y,k^*_{s,y} , \nabla_x J(s,y)).
\end{align}

Moreover, since $V^*$ is the viscosity solution to the HJB PDE by Lemma~\ref{lem: unique solution to HJB PDE}, by definition we have that
\vspace{-0.2cm} \begin{align} \label{12}
\frac{\partial}{\partial s} V^*(s,y)+  H(s,y,k^*_{s,y} , \nabla_x V^*(s,y))=0.
\end{align}

Combining Eqs.~\eqref{11} and \eqref{12}, {applying H\"older's inequality ($a^\top b \le \lVert a \rVert_{\infty} \lVert b \rVert_{1} \le \lVert a \rVert_{2} \lVert b \rVert_{1}$), using the bound $\lVert f(t,x,u) \rVert_2 \le \beta_f (1+\lVert x \rVert_2)$} (Defn.~\ref{ass: conditions on OCP}), and the fact $x_{u,x_0}(t) \in B_R(0)$ (by Gronwall Bellman), it follows that
\begin{align} \nonumber
& E_J(s,y)  \le  \frac{\partial}{\partial s} J(s,y) - \frac{\partial}{\partial s} V^*(s,y)  \\ \nonumber
& \quad \qquad \qquad  + ( \nabla_x J(s,y) - \nabla_x V^*(s,y))^{\top} f(s,y,k^*_{s,y})\\ \nonumber
& \le \left|\frac{\partial}{\partial s} J(s,y) - \frac{\partial}{\partial s} V^*(s,y)\right| \\  \nonumber
& + \hspace{-5mm}\sup_{(t,x,u) \in [0,T] \times B_R(0) \times U} \hspace{-0.8cm} \|f(t,x,u)\|_2  \sum_{i=1}^n \bigg| \frac{\partial }{\partial y_i} (J(s,y)  -V^*(s,y) ) \bigg|  \\  \nonumber
& \le \max \bigg\{ 1, \beta_f (1+R) \bigg\} \bigg (\left|\frac{\partial}{\partial s} J(s,y) - \frac{\partial}{\partial s} V^*(s,y)\right| \\ \label{42}
& \qquad + \sum_{i=1}^n \bigg| \frac{\partial }{\partial y_i} (J(s,y)  -V^*(s,y) ) \bigg| \bigg ).
\end{align}


\underline{\textbf{Part 2 of Proof:}} If ($s,y) \in S_J \cap S_{V^*}$ and {$\pi_J$} satisfies Eq.~\eqref{eq: controller synth}, then {by Eq.~\eqref{pfeq: simplfication of E ae}}
\vspace{-0.2cm} \begin{align} \label{21}
E_J(s,y) = \frac{\partial}{\partial s}J(s,y) & + H(s,y,\pi_J(s,y), \nabla_x J(s,y)).  
\end{align}
Moreover, since $V^*$ is the viscosity solution to the HJB PDE~\eqref{eqn: general HJB PDE}, by Lemma~\ref{lem: unique solution to HJB PDE}, we have that
\vspace{-0.2cm} \begin{align} \label{22}
&	\frac{\partial}{\partial s} V^*(s,y) +	H(s,y,\pi_J(s,y), \nabla_x V^*(s,y))\\ \nonumber
& \ge  	\frac{\partial}{\partial s} V^*(s,y) +	\inf_{u \in U} H(s,y,u, \nabla_x V^*(s,y)) =  0. 
\end{align}
Combining Eqs~\eqref{21} and \eqref{22} it follows by a similar argument to Eq.~\eqref{42} that,
\ifarxiv \begin{align} \nonumber
& E_J(s,y)  \ge     \frac{\partial}{\partial s} J(s,y) - \frac{\partial}{\partial s} V^*(s,y)  \\ \nonumber
& \quad + ( \nabla_x J(s,y) - \nabla_x V^*(s,y))^{\top} f(s,y,\pi_J(s,y))\\ \nonumber
& \ge -|\frac{\partial}{\partial s} J(s,y) - \frac{\partial}{\partial s} V^*(s,y)| \\ \nonumber
& - \hspace{-5mm}\sup_{(t,x,u) \in [0,T] \times B_R(0) \times U} \hspace{-0.8cm} \|f(t,x,u)\|_2\sum_{i=1}^n \bigg| \frac{\partial }{\partial y_i} (J(s,y)  -V^*(s,y) ) \bigg|\\ \nonumber
& \ge -\max \bigg\{ 1, \beta_f (1+R) \bigg\} \bigg (\left|\frac{\partial}{\partial s} J(s,y) - \frac{\partial}{\partial s} V^*(s,y)\right| \\ \label{41}
& \qquad + \sum_{i=1}^n \bigg| \frac{\partial }{\partial y_i} (J(s,y)  -V^*(s,y) ) \bigg| \bigg ).
\end{align}
\else
\begin{align} \nonumber
& E_J(s,y)  \\ \nonumber
& \ge -\max \bigg\{ 1, \beta_f (1+R) \bigg\} \bigg (\left|\frac{\partial}{\partial s} J(s,y) - \frac{\partial}{\partial s} V^*(s,y)\right| \\ \label{41}
& \qquad + \sum_{i=1}^n \bigg| \frac{\partial }{\partial y_i} (J(s,y)  -V^*(s,y) ) \bigg| \bigg ).
\end{align} \fi
\vspace{-0.2cm}

\underline{\textbf{Part 3 of Proof:}}
{From Eqs~\eqref{42} and~\eqref{41}} it follows for any $(s,y) \in S_J \cap S_{V^*}$ we have 
\begin{align*} \nonumber
   & |E_J(s,y)| \le \max \bigg\{ 1, \beta_f (1+R) \bigg\} \bigg (\left|\frac{\partial}{\partial s} (J(s,y) -  V^*(s,y))\right| \\  & \hspace{2cm} + \sum_{i=1}^n \bigg| \frac{\partial }{\partial y_i} (J(s,y)  -V^*(s,y) ) \bigg| \bigg ).
\end{align*}
and hence since $([0,T] \times B_R(0))/(S_J \cap S_{V^*})$ has zero measure it follows that,
\vspace{-0.4cm}\begin{align} \nonumber
    & \esssup_{(s,y) \in [0,T] \times B_R(0)} \{ |E_J(s,y)| \}\le \max \bigg\{ 1, \beta_f (1+R) \bigg\} \\ \nonumber & \qquad \times  \esssup_{(s,y) \in [0,T] \times B_R(0)} \bigg(\left|\frac{\partial}{\partial s} (J(s,y) -  V^*(s,y))\right| \\ \label{pfeq: ess sup} & \hspace{2cm} + \sum_{i=1}^n \bigg| \frac{\partial }{\partial y_i} (J(s,y)  -V^*(s,y) ) \bigg| \bigg ).
\end{align}

{Combining Inequalities~\eqref{pfeq: loss bound} and \eqref{pfeq: ess sup} we obtain~\eqref{ineq:loss functin}, completing the proof.
\vspace{-0.2cm}\begin{align*}
&L( u_{\pi_J,x_0},x_0)  
\le  2T\max \bigg\{ 1, \beta_f (1+R) \bigg\}\\ \nonumber & \qquad \times \esssup_{(s,y) \in [0,T] \times B_R(0)}   \bigg (\left|\frac{\partial}{\partial s} (J(s,y) -  V^*(s,y))\right| \\  & \hspace{2cm} + \sum_{i=1}^n \bigg| \frac{\partial }{\partial y_i} (J(s,y)  -V^*(s,y) ) \bigg| \bigg )\\
&+  2 \hspace{-0.05cm}  \sup_{y \in B_R(0)}\{ |V^*(T,y)- J(T,y)| \}\\
& \le C  \|J- V^*\|_{W^{1,\infty}([0,T] \times B_R(0))},
\end{align*}
where $C:= 2\max \left\{1,T, T {\beta_f(1+R)} \right\}$.}
\end{proof}

\vspace{-0.5cm} 
{\section{Interpretation of Performance Bounds}

Prop.~\ref{prop: Loss bounded by HJB res} shows that direct minimization of the HJB residual is ideal for controller synthesis. In practice, however, directly minimizing this residual is often computationally intractable due to the residual itself containing an implicit inner optimization problem ($\inf_{u \in U}$). Consequently, many numerical methods must construct candidate VFs via other means. Theorem~\ref{thm: performance bounds} bridges this gap: if a numerical scheme approximates the true VF, securing a small error in the $W^{1,\infty}$ norm provides a guarantee that the resulting controller will perform well. Later, we will see in Example~\ref{ex:lib} that a small error in the weaker $W^{1,p}$ norm for $p<\infty$ alone (without HJB residual minimization) is not strong enough to produce good controller performance.


Importantly, while a small $W^{1,\infty}$ error is a sufficient condition for good controller performance, Example~\ref{ex:ortho_gradient}, presented next, demonstrates that it is not strictly necessary. It is possible to construct a candidate VF with an arbitrarily large $W^{1,\infty}$ error that still yields a arbitrarily near optimal controller. This occurs when the gradient approximation error is orthogonal to the system's vector field, effectively hiding the error from the Hamiltonian. In such cases, the HJB residual remains zero despite the divergent gradient. However, generic numerical approximation schemes do not generally restrict gradient errors exclusively to directions orthogonal to the system's vector field. Therefore, Theorem~\ref{thm: performance bounds} provides a universal performance guarantee, bounding controller performance by the Sobolev VF error.




\begin{ex} \label{ex:ortho_gradient}
Consider an OCP with tuple $\{c,g,f,U,T\}$ defined by $c \equiv 0$, $g(x)=\|x\|_2^2$, $U=[-1,1]$, $T=1$ and $f(x,u)=[1, -1]^\top u \|x\|_2^2$. Let $V^*$ be the true VF and for any $\epsilon >0$, consider the admissible candidate value function $J_\eps(t,x):=V^*(t,x) + \sqrt{\eps} \sin((x_1+x_2)/ \eps)$. 
First, we observe that $\nabla_x J_\eps(t,x) = \nabla_x V^*(t,x) + \eps^{-1/2} \cos((x_1+x_2)/\eps) [1,1]^\top$, so $\lim_{\epsilon\rightarrow 0}\norm{\nabla_x J_\eps-\nabla_x V^*(t,x)}=\infty$ -- implying the error in the $W_{1,\infty}$ norm becomes arbitrarily large. 

Despite this divergence, the HJB residual meanwhile becomes arbitrarily small. To see this, we note that 
$\frac{\partial}{ \partial t} J_\eps(t,x)= \frac{\partial}{ \partial t}V^*(t,x)$ and $\nabla_x J_\eps^\top f = \nabla_x V^{*\top} f + \eps^{-1/2}\cos((x_1+x_2)/\eps) [1,1] [1, -1]^\top u \|x\|_2^2 = \nabla_x V^{*\top} f$.  We conclude that 
\begin{align*}
&\left|\frac{\partial}{\partial t}J_{\epsilon}(t,x)  + \inf_{u \in U} \{H(t,x,u, \nabla_x J_{\epsilon}(t,x)) \}\right|\\
&=\left|\frac{\partial}{\partial t}V^*(t,x)  + \inf_{u \in U} \{H(t,x,u, \nabla_x V^*(t,x)) \}\right|=0.
\end{align*}
Furthermore,  $\lim_{\epsilon \rightarrow 0}\| J_\eps(t,\cdot) - V^*(t,\cdot)\|_{L^\infty([0,T],\R^n)} = 0$, implies $\lim_{\eps \to 0} J_\eps(T,x) = V^*(T,x) = g(x)$. Applying Prop.~\ref{prop: Loss bounded by HJB res}, we conclude that $\lim_{\epsilon \rightarrow 0}L( u_{\pi_{J_\epsilon},x_0},x_0)=0$ -- implying that the performance of the sequence of resulting controllers converges to optimality.
%
%

\end{ex}

Finally, we note that the construction of the bounds in Prop.~\ref{prop: Loss bounded by HJB res} and Thm.~\ref{thm: performance bounds} use the Gronwall Bellman inequality to bound growth of the solution -- resulting in an exponential dependency in time-horizon, $T$. Such exponentials may be avoided if a bound on the forward reachable set is known a priori.} Furthermore, because the suboptimality of {any} controller cannot exceed the gap between the worst-case (maximum cost) VF, we denote as $\mathcal{V}^*$, and the optimal VF, $V^*$, we can refine Eq.~\eqref{ineq:loss functin}:
\vspace{-0.1cm}\begin{align} \nonumber
L(u_{\pi_J,x_0},x_0) \leq \min \bigg\{ C & \| J - V^* \|_{W^{1,\infty}([0,T] \times B_R(0))},\\ \label{eq: Improved per bound}
&  \mcl V^* (0,x_0) - V^*(0,x_0) \bigg\}.
\end{align}

\vspace{-0.4cm}
\section{Numerical Examples} \label{sec: numerical ex}
\ifarxiv We proceed by synthesizing several closed-loop controllers from various candidate VFs and assess their performance at minimizing a given cost function, both numerically and in terms of the theoretical performance bounds derived in Eq.~\eqref{eq: Improved per bound}. Trajectories are computed using MATLAB's \texttt{ode45} function. \fi


\begin{figure*}[h!]
\begin{subfigure}[t]{ 0.333 \textwidth} 
\begin{tikzpicture}[>=Stealth, scale=2.5]

\definecolor{trueVblue}{RGB}{0,114,189}
\definecolor{trapVred}{RGB}{217,83,25}

\pgfmathsetmacro{\tval}{0}
\pgfmathsetmacro{\eps}{0.1}
\pgfmathsetmacro{\deltaVal}{0.1}
\pgfmathsetmacro{\xc}{0.75}
\pgfmathsetmacro{\amp}{pow(\eps,\deltaVal)}

\pgfmathsetmacro{\yBase}{\xc*exp(-1)}
\pgfmathsetmacro{\yTrap}{\yBase-\amp*exp(-1)}

\draw[->, thick] (-0.5,0) -- (1.5,0) node[below] {$x$};
\draw[->, thick] (0,-0.8) -- (0,1) node[above] {};

\draw[
dashed,
very thick,
color=trueVblue,
domain=-0.3:1.5,
samples=20
]
plot (\x,{
(\x>=0)*(\x*exp(-1)) + (\x<0)*(\x*exp(1))
});

\draw[
dotted,
very thick,
color=trapVred,
domain=-0.3:1.5,
samples=80
]
plot (\x,{
(\x>=0)*(\x*exp(-1)) + (\x<0)*(\x*exp(1))
- \amp *
(abs((\x-\xc)/\eps) < 1) *
exp( ( (abs((\x-\xc)/\eps) < 1)/(pow((\x-\xc)/\eps,2) - 1 ) ) )	});

\begin{scope}[shift={(0.7,0.7)}]
\draw[thick, fill=white, fill opacity=0.9, draw opacity=1] 
(-0.05,-0.05) rectangle (0.8,0.35);

\draw[dashed, very thick, color=trueVblue] 
(0,0.2) -- (0.2,0.2);
\node[right, color=trueVblue] at (0.22,0.2) {$V^*(0,x)$};

\draw[dotted, very thick, color=trapVred] 
(0,0.05) -- (0.2,0.05);
\node[right, color=trapVred] at (0.22,0.05) {$\hat{V}_1(0,x)$};
\end{scope}

\coordinate (Origin) at (0,0);
\fill (Origin) circle (0.8pt);

\draw[<-]
($(Origin)+(-0.1,0.02)$) -- ++(-0.3,0.5)
node[above, align=center, font=\footnotesize]
{Non-diff.\\region ($x=0$)};

\draw[<-]
(0.6,0.31) -- (0.3,0.6)
node[above, align=center, font=\scriptsize, color=trueVblue]
{$\nabla V^* > 0$\\$\pi_{V^*} =-1$};

\coordinate (TrapSlope) at (0.65,0.1);
\draw[<-]
(TrapSlope) -- (0.5,-0.2)
node[below, align=center, font=\scriptsize, color=trapVred]
{$\nabla \hat{V}_1 < 0$\\$\pi_{\hat{V}_1} = +1$};

\draw[
<->, 
thin
]
(\xc-\eps,\yBase+0.1) -- (\xc+\eps,\yBase+0.1)
node[midway, above=3pt, font=\tiny]
{Width=$\epsilon$};


\draw[<->, thin]
(\xc+1.2*\eps,\yBase) -- (\xc+1.2*\eps,\yTrap)
node[midway, right, font=\tiny]
{Height$ \sim \epsilon^\delta$};

\end{tikzpicture}
\vspace{-20pt}
\caption{}
\label{fig: Illistrative example}
\end{subfigure}%
\begin{subfigure}[t]{0.3333 \textwidth}
\includegraphics[width=\linewidth, trim = {3cm 8cm 2.5cm 9cm}, clip]{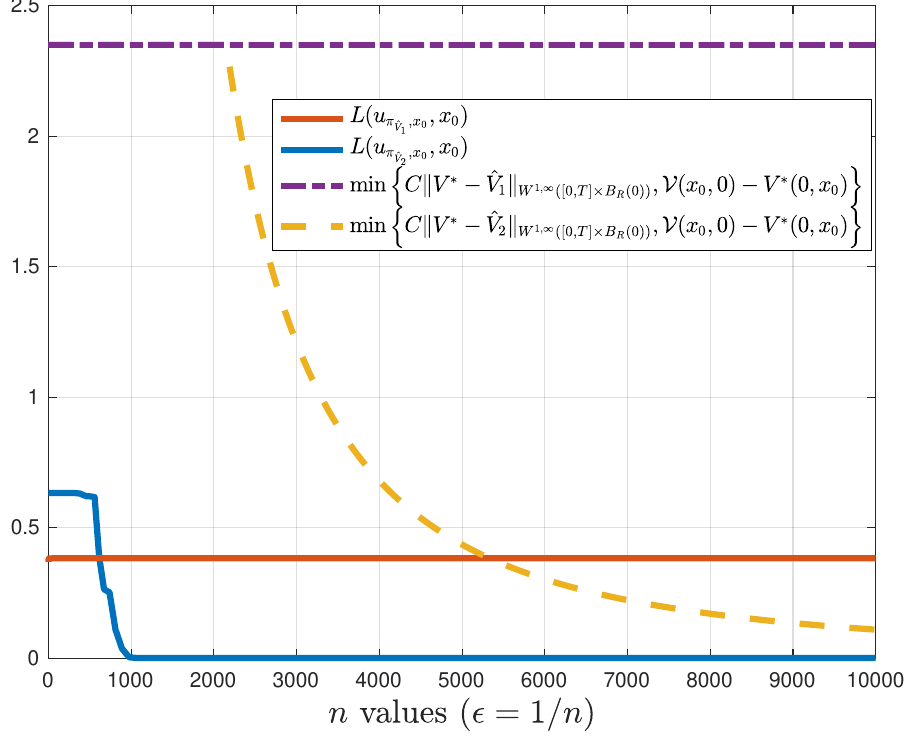}
\vspace{-20pt}
\caption{}
\label{fig: lokta}
\end{subfigure}%
\begin{subfigure}[t]{ 0.333 \textwidth}
\includegraphics[width=\linewidth, trim =  {1.5cm 0cm 1.5cm 1cm}, clip]{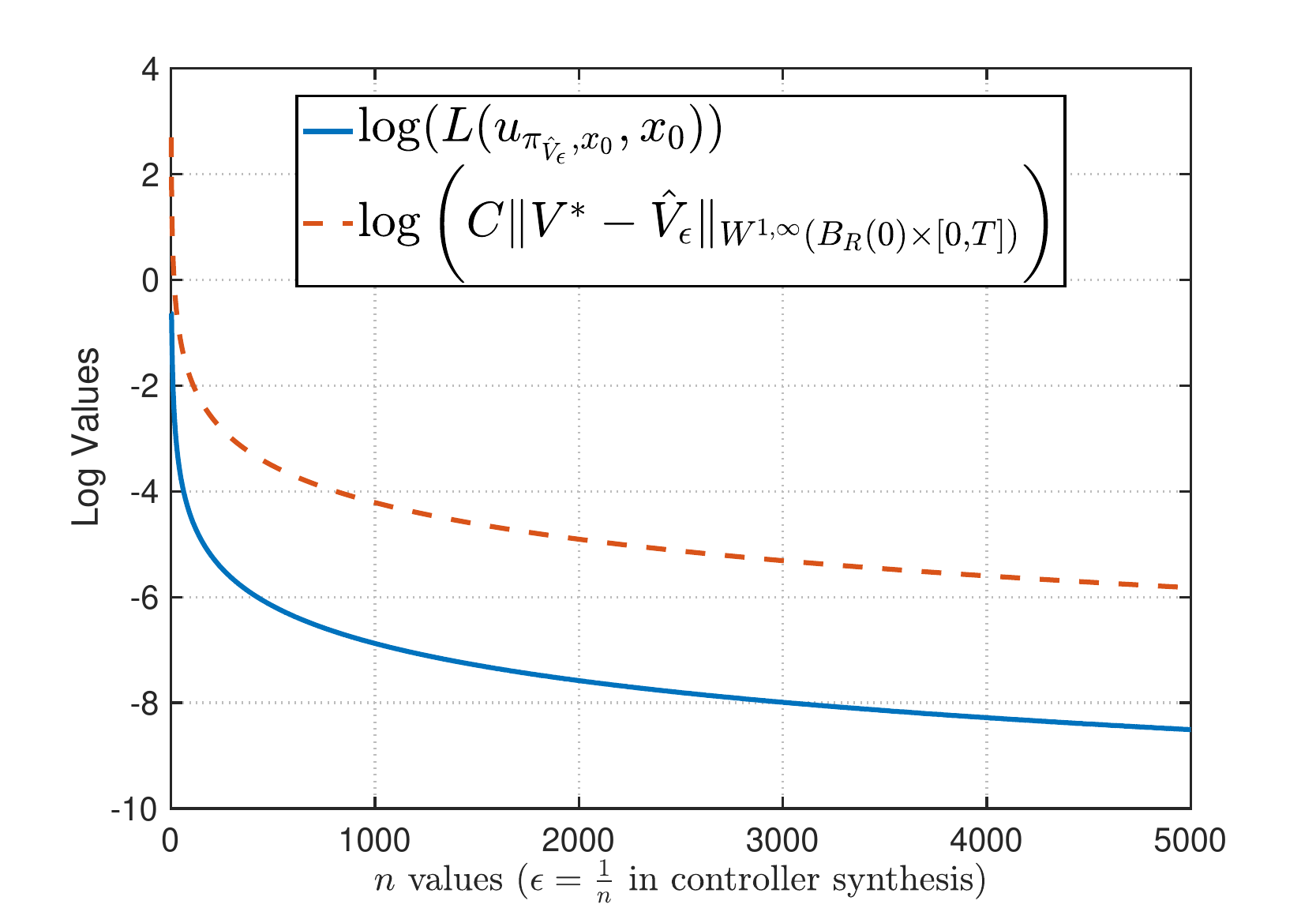}
\vspace{-20pt}
\caption{}
\label{fig: coupled linear ODE}
\end{subfigure}%
\vspace{-10pt}
\caption{\footnotesize (\protect\subref{fig: Illistrative example}) {Graph showing plots of $V^*$  and $\hat{V}_1$ from Eqs~\eqref{analytical soln} and \eqref{eq: V1 and V2} respectively.}
(\protect\subref{fig: lokta})	Graph showing numerically calculated performance gaps, $L(u_{\pi_{\hat{V}_1},x_0},x_0)$ and $L(u_{\pi_{\hat{V}_2},x_0},x_0)$, associated with the approximate VFs given in Eq.~\eqref{eq: V1 and V2} for decreasing $\eps>0$, along with the theoretical bounds for these performance gaps given in Eqs~\eqref{eq: ex1 theoretical performance bounds} and~\eqref{eq: ex1 theoretical performance bounds V 2}.
(\protect\subref{fig: coupled linear ODE}) 	Numerically calculated logarithms of performance gap, $L(u_{\pi_{\hat{V}_\eps},x_0},x_0)$, associated with the approximate VFs given in Eq.~\eqref{eq: approx VF example 2} for decreasing $\eps>0$, along with the logarithm of the theoretical performance gap bound given in Eq.~\eqref{eq: perf bound ex 2}.}
\vspace{-10pt}
\end{figure*}

\begin{ex} \label{ex:lib}
Let us consider the OCP associated with the tuple $\{c,g,f,U,T\}$, where $c(t,x,u) \equiv 0$, $g(x)=x$, $f(t,x,u)=xu$, $U=[-1,1]$, and $T=1$. It was shown in \cite{liberzon2011calculus} that the VF associated with this tuple can be found analytically:
\vspace{-0.2cm}
\begin{equation} \label{analytical soln}
V^*(t,x)= \begin{cases} 
\exp(t-1)x & \text{if } x>0,\\
\exp(1-t)x & \text{if } x<0,\\
0 & \text{if } x=0. 
\end{cases}
\end{equation}
\ifarxiv By a similar argument, we can find the VF that maximizes cost:
\vspace{-0.5cm} 
\begin{equation} \label{analytical soln 2}
\hspace{1.5cm} \mcl V^*(t,x)= \begin{cases} 
\exp(1-t)x & \text{if } x>0,\\
\exp(t-1)x & \text{if } x<0,\\
0 & \text{if } x=0. 
\end{cases}
\end{equation} \fi

Let us now consider the following candidate VFs:
\begin{align} \label{eq: V1 and V2}
\hat{V}_1(t,x) & = V^*(t,x) -{100 \eps^\delta \eta ((x-0.75)/\eps)}, \\ \nonumber
\hat{V}_2(t,x) & = V^*(t,x) + \eps^2 \sin(10^6x),
\end{align}
{where $\eta(x):= \begin{cases}
\exp \left( \frac{1}{x^2 -1} \right) \text{ if } x \in (-1,1)\\
0 \text{ otherwise.}
\end{cases}$ is the standard bump function from mollification theory and $\delta \in (0,1)$. It is well known that $\eta \in C^\infty(\R, \R)$~\cite{evans2010partial} with compact support; hence $C_{0,p}:=\| \eta \|_{L^p(\R)}< \infty$ and $C_{1,p}:=\| \eta' \|_{L^p(\R)}< \infty$ for all $p \ge 1$, including $p=\infty$. \ifarxiv For $p=\infty$ these constants can be numerically approximated using unconstrained nonlinear optimization.

To bound the error between $V^*$ and $\hat{V}_1$, we evaluate the Sobolev norm of the difference over the spacetime domain $[0,T] \times \R$. Let the difference be denoted by $\mathcal{E}_1(x) := V^*(t,x) - \hat{V}_1(t,x) = 100 \eps^\delta \eta \left( \frac{x-0.75}{\eps} \right)$. 

For $1 \le p < \infty$, and noting $\mathcal{E}_1$ is independent of time, it follows,}
\begin{align*}
   & \| \mathcal{E}_1 \|_{W^{1,p}([0,T] \times \R)}^p  \\
    &= \int_0^T \int_{\R} |\mathcal{E}_1 (x)|^p \,dx\,dt + \int_0^T \int_{\R} | \mathcal{E}_1 '(x)|^p \,dx\,dt \\
    &= T \int_{\R} \left| 100 \eps^\delta \eta\left(\frac{x-0.75}{\eps}\right) \right|^p dx  \\
    & \qquad + T \int_{\R} \left| 100 \eps^{\delta-1} \eta'\left(\frac{x-0.75}{\eps}\right) \right|^p dx \\
    &=  100^p T\eps^{\delta p} \int_{\R} |\eta(y)|^p \eps \, dy +  100^p T \eps^{(\delta - 1)p} \int_{\R} |\eta'(y)|^p \eps \, dy \\
    &=  100^p T\left( \eps^{\delta p + 1} C_{0,p}^p + \eps^{(\delta - 1)p + 1} C_{1,p}^p \right).
\end{align*}
Taking the $p$-th root gives,
\begin{align*}
   & \| \mathcal{E}_1 \|_{W^{1,p}([0,T] \times \R)} \\
    & \qquad = 100 T^{\frac{1}{p}} \left( \left(\eps^{\delta + \frac{1}{p}} C_{0,p}\right)^p + \left(\eps^{\delta - 1 + \frac{1}{p}} C_{1,p}\right)^p \right)^{\frac{1}{p}}.
\end{align*}
Using the inequality $(a^p + b^p)^{\frac{1}{p}} \le a + b$ for $a,b \ge 0$ and $p \ge 1$, we obtain the upper bound:
\begin{equation} \label{eq: V star - V}
    \| V^* - \hat{V}_1 \|_{W^{1,p}([0,T] \times \R)} \le 100 T^{\frac{1}{p}} \left( \eps^{\delta + \frac{1}{p}} C_{0,p} + \eps^{\delta - 1 + \frac{1}{p}} C_{1,p} \right).
\end{equation}

For the case where $p = \infty$, 
\begin{align*}
  &  \| \mathcal{E}_1 \|_{W^{1,\infty}([0,T] \times \R)} \\
    &= \operatorname*{ess\,sup}_{(t,x) \in [0,T] \times \R} |\mcl E_1 (x)|  + \operatorname*{ess\,sup}_{(t,x) \in [0,T] \times \R} |\mcl E_1'(x)| \\
    &= 100 \eps^\delta \operatorname*{ess\,sup}_{x \in \R} \left| \eta\left(\frac{x-0.75}{\eps}\right) \right| \\
    & \qquad \qquad + 100 \eps^{\delta-1} \operatorname*{ess\,sup}_{x \in \R} \left| \eta'\left(\frac{x-0.75}{\eps}\right) \right| \\
    &= 100 \eps^\delta \|\eta\|_{L^\infty(\R)} + 100 \eps^{\delta - 1} \|\eta'\|_{L^\infty(\R)} \\
    &= 100 \left( \eps^\delta C_{0,\infty} + \eps^{\delta - 1} C_{1,\infty} \right).
\end{align*}

We next bound the difference between $V^*$ and $\hat{V}_2$. Denote $\mcl E_2(x):=V^*(t,x)-\hat{V}_2(t,x)=\eps^2 \sin(10^6 x)$. It follows that,
\begin{align*}
      &  \| \mathcal{E}_2 \|_{W^{1,\infty}([0,T] \times \R)} \\
    &= \operatorname*{ess\,sup}_{(t,x) \in [0,T] \times \R} |\mcl E_2 (x)|  + \operatorname*{ess\,sup}_{(t,x) \in [0,T] \times \R} |\mcl E_2'(x)| \\
    & = \eps^2 + 10^6 \eps^2.
\end{align*}
Applying the performance bound in Eq.~\eqref{eq: Improved per bound} for $x_0=1$, $T=1$ and $\beta_f=1$, we deduce:
\vspace{-0.2cm}
\begin{align} \nonumber
& L(u_{\hat{V}_1},x_0) \le \min \bigg\{ 200 \max\{1,T,T\beta_f(1+|x_0|)e^{\beta_f T} \}\\ \nonumber
& \hspace{2cm} \times (\eps^\delta C_{0,\infty} +\eps^{\delta -1} C_{1,\infty} ), e-e^{-1} \bigg\} \\ \label{eq: ex1 theoretical performance bounds}
& =\min \left\{400 e \left( C_{0,\infty} \eps^\delta + \frac{C_{1,\infty}}{{\eps^{1-\delta}}} \right), (e-e^{-1}) \right\},\\ \nonumber
& L(u_{\hat{V}_2},x_0)  \le \min \bigg\{ 2 \max\{1,T,T\beta_f(1+|x_0|)e^{\beta_f T} \}\\ \nonumber
& \hspace{2cm} \times (1+10^6) \eps^2, (e-e^{-1})  \bigg\} \\  \label{eq: ex1 theoretical performance bounds V 2} 
& = \min \left\{4e(1+10^6) \eps^2 , (e-e^{-1}) \right\}.
\end{align}

\else
Therefore\footnote{See full derivation details on Arxiv~\cite{jones2025bounding}}:}
\vspace{-0.1cm}
\begin{align} \nonumber
&	\| V^* \hspace{-0.05cm}  - \hspace{-0.05cm}  \hat{V}_1 \|_{W^{1,p}([0,T] \times \R)} \hspace{-0.1cm} {\le \hspace{-0.05cm} 100 (\eps ^{\delta+1/p} C_{0,p} \hspace{-0.05cm}+ \hspace{-0.05cm} \eps^{\delta-1+1/p} C_{1,p})} ,\\ \label{eq: V star - V}
&	\| V^*-\hat{V}_2 \|_{W^{1,\infty}([0,T] \times \R)}  =(1+10^6) \eps^2.
\end{align}
Applying the performance bound in Eq.~\eqref{eq: Improved per bound} for $x_0=1$ and $\beta_f=1$, we deduce:
\vspace{-0.2cm}
\begin{align} \label{eq: ex1 theoretical performance bounds}
L(u_{\pi_{\hat{V}_1},x_0},x_0)  \hspace{-0.05cm} & {\le \hspace{-0.05cm}\min \left\{ \hspace{-0.05cm} 400 e \left( \hspace{-0.05cm} C_{0,\infty} \eps^\delta \hspace{-0.05cm}+ \hspace{-0.05cm}\frac{C_{1,\infty}}{{\eps^{1-\delta}}} \hspace{-0.05cm} \right) \hspace{-0.05cm}, (e\hspace{-0.05cm} - \hspace{-0.05cm} e^{-1}) \hspace{-0.05cm} \right\}},\\  \label{eq: ex1 theoretical performance bounds V 2} 
L(u_{\pi_{\hat{V}_2},x_0},x_0) & \le \min \left\{4e(1+10^6) \eps^2 , (e-e^{-1}) \right\}.
\end{align}  \fi
{Observe via Eq.~\eqref{eq: V star - V} that for any finite $p < \infty$, we can select $\delta \in (1 - 1/p, 1)$ such that $\hat{V}_1 \to V^*$ in the $W^{1,p}$ norm. Despite this, the synthesized controller fails to yield trajectories that converge to the optimal cost. This discrepancy arises because $\hat{V}_1$ does not converge to $V^*$ in the stronger $W^{1,\infty}$ norm. As illustrated in Fig.~\ref{fig: Illistrative example}, for every $\eps>0$ and $\delta \in (0,1)$, the bump function perturbation creates a local minimum in $\hat{V}_1$ that acts as a ``trapping region.'' Consequently, the induced controller stabilizes the trajectory around the center of the bump ($x=0.75$) rather than optimally driving the state to minimize the terminal cost $g(x)=x$.

This phenomenon is confirmed numerically in Fig.~\ref{fig: lokta} for $\delta=0.9$, where the cost does not converge to zero as $\eps \to 0$. The right-hand side of the bound in Eq.~\eqref{eq: ex1 theoretical performance bounds} remains constant at $e-e^{-1}\approx 2.35$ as $\eps \to 0$. While this provides a valid upper bound, it indicates that the controller is not guaranteed to achieve optimal performance.}

Conversely, Eq.~\eqref{eq: ex1 theoretical performance bounds V 2} demonstrates that the controller synthesized by $\hat{V}_2$ will converge to the optimal cost, since $L(u_{\pi_{\hat{V}_2},x_0},x_0) \to 0$ as $\eps \to 0$. In Fig.~\ref{fig: lokta}, we observe that $u_{\pi_{\hat{V}_2},x_0}$ numerically converges to optimal performance at a faster rate than the theoretical error bound derived in Eq.~\eqref{eq: ex1 theoretical performance bounds V 2}.
\end{ex}
\begin{ex}
Let us consider the OCP associated with the following tuple $\{c,g,f,U,T\}$, where $c(t,x,u) = (1+2(T-t))x^2 + 2(T-t)|x|$, $g(x)=0$, $f(t,x,u)=-x +u $, $U=[-1,1]$, and $T=1$. Since, $\inf_{u \in [-1,1]} \{\nabla_x V(t,x) u \}=-|\nabla_x V(t,x)|$ for any differentiable $V$ it follows that the associated HJB PDE is
\vspace{-0.1cm}\begin{align} \label{eq: HJB for ex}
\frac{\partial}{\partial t} V(t,x)+& (1+2(T-t))x^2 + 2(T-t)|x| \\ \nonumber
& - \nabla_x V(t,x) x - |\nabla_x V(t,x)|=0.
\end{align}
We can then deduce $V(t,x)=x^2(T-t)$ is the VF for this OCP by substitution into Eq.~\eqref{eq: HJB for ex}. Hence the optimal controller is given by
\vspace{-0.5cm}\begin{align}
\pi^*(t,x)= -\sign(2x(T-t))=\begin{cases}
-1 \text{ if } x>0,\\
0 \text{ if } x=0,\\
1 \text{ if } x<0.
\end{cases}
\end{align}
Let us next consider a sequence of approximate VFs
\vspace{-0.1cm}\begin{align} \label{eq: approx VF example 2}
\hat{V}_\eps (t,x) & = x^2(T-t) + \eps x.
\end{align}
The candidate VF, $	\hat{V}_\eps$, yields the following controller
\vspace{-0.1cm}\begin{align}
\hat{\pi}_{\hat{V}_\eps}\hspace{-0.05cm} (t, \hspace{-0.05cm}x)\hspace{-0.1cm}=\hspace{-0.05cm} -\hspace{-0.05cm} \sign(2x(T \hspace{-0.1cm}- \hspace{-0.05cm}t) \hspace{-0.05cm} + \hspace{-0.05cm} \eps )\hspace{-0.1cm}= \hspace{-0.05cm}\begin{cases}
-1 \hspace{-0.05cm} \text{ if } \hspace{-0.05cm} x\hspace{-0.05cm} >\hspace{-0.05cm} \frac{-\eps}{2(T-t)},\\
0 \text{ if } x=\frac{-\eps}{2(T-t)},\\
1 \text{ if } x<\frac{-\eps}{2(T-t)}.
\end{cases}
\end{align}
For $x_0=0$ we have that $R:=(1+\|x_0\|_2)e^{\beta_f T}-1=e-1$. Hence, $	\| V^*-\hat{V}_\eps \|_{W^{1,\infty}([0,T] \times B_R(0) )} \le (R+1)\eps \to 0$ but $	\lim_{\eps \to 0}\| \pi^*-\hat{\pi}_{\hat{V}_\eps}\|_{L^{\infty}([0,T] \times B_R(0), \R)} \ne 0$ since $\pi^*(t,0)=0$ but $\hat{\pi}_{\hat{V}_\eps}(t,0)=-1$ for all $\eps>0$. Despite the fact that the controller $\hat{\pi}_{\hat{V}_\eps}$ does not converge to the optimal controller uniformly, Theorem~\ref{thm: performance bounds} still shows that the associated performance/cost of applying the controller converges to zero. This is shown numerically in Fig.~\ref{fig: coupled linear ODE}. The figure also shows the theoretical performance bound found by substituting $R=e-1$, $\beta_f=1$ into Eq.~\eqref{ineq:loss functin} to get
\vspace{-0.6cm} \begin{flalign} \label{eq: perf bound ex 2}
&\text{ } \qquad \qquad\qquad \qquad  L( u _{\pi_{\hat{V}_\eps},x_0},x_0) \le 2 e^2 \eps. 
\end{flalign}
\end{ex}
\vspace{-0.35cm}
\section{Conclusion}
\vspace{-0.2cm}
We have shown that when approximate value functions are used to synthesize controllers, the suboptimality of the resulting controller is bounded by {the $L^\infty$ norm of the HJB residual, which is, in turn, bounded by the Sobolev $W^{1,\infty}$ norm of the error between the candidate value function and true value function. Consequently, numerical methods for solving the HJB equation which converge in the $W^{1,\infty}$ norm will yield controllers which minimize the associated cost arbitrarily close to the true optimum. In contrast, we also showed that convergence in $W^{1,p}$ for finite $p$ does not necessarily guarantee optimal performance. These findings suggest that future work should prioritize minimizing the $W^{1,\infty}$ distance to the true value function or $L^\infty$ HJB residual error over other metrics and objectives.}

\bibliographystyle{ieeetr}
\bibliography{bibliography}

\ifarxiv 
\section{Appendix A}
\begin{lemma} \label{lem: ess sup bound}
		Let $E : [0,T] \times \mathbb{R}^n \to \mathbb{R}$ be a lower semicontinuous function. For any measurable function $x: [0,T] \to B_R(0) \subset \mathbb{R}^n$, for some $R>0$, the following inequality holds:
		\begin{equation}
			\operatorname*{ess\,sup}_{t \in [0,T]} E(t, x(t)) \le \operatorname*{ess\,sup}_{(t,y) \in [0,T] \times B_R(0)} E(t, y).
		\end{equation}
	\end{lemma}
	
	\begin{proof}
		Let $M = \operatorname*{ess\,sup}_{(t,y) \in [0,T] \times B_R(0)} E(t, y)$. Suppose, for the sake of contradiction, that the inequality does not hold. Then, there exists $\varepsilon > 0$ such that the set
		\[
		S = \{ t \in [0,T] \mid E(t, x(t)) > M + \varepsilon \}
		\]
		has positive measure, i.e., $\mu(S) > 0$. Hence, there exists $t_0 \in S \cap (0,T)$. Consider the point $p_0 = (t_0, x(t_0)) \in (0,T) \times B_R(0)$. Since $E$ is lower semicontinuous, the superlevel set
		\[
		V = \{ (t, y) \in (0,T) \times B_R(0) \mid E(t, y) > M + \varepsilon \}
		\]
		is open. Since $p_0 \in V$, $V$ is a non-empty open set. In Euclidean spaces equipped with the Lebesgue measure, every non-empty open set has positive measure. 
		
		This implies that $E(t, y) > M + \varepsilon$ on a set of positive measure in the product space, which contradicts the definition of $M$ as the essential supremum. Thus, the assumption must be false, and the inequality holds.
\end{proof}

\section{Appendix B: A Generalized Chain Rule}

In this appendix we provide a version of the chain rule for the composition of Lipschitz functions and absolutely continuous trajectories. We discuss the intuition behind this chain rule generalization and cases in which the rule cannot be applied.
\begin{lemma} \label{lem:chain_rule}
Consider functions $J \in \text{Lip}([0,T] \times \mathbb{R}^n , \mathbb{R})$ and $x \in {\text{AC}}([0,T], \mathbb{R}^n)$ satisfying $\dot{x}(t) = f(t,x(t))$ almost everywhere on $[0,T]$.
Then, for almost every $t \in [0,T]$ and any $(p,q) \in D^\pm J(t,x(t))$, we have
\begin{align} \label{eq: chain rule}
\frac{d}{dt} J(t,x(t)) &= p + q^\top f(t,x(t)).
\end{align}
\end{lemma}

\begin{proof}
Since $x(\cdot)$ is absolutely continuous and $J$ is locally Lipschitz, the composite function $g(t) := J(t,x(t))$ is absolutely continuous on $[0,T]$. Consequently, the classical derivative $g'(t)$ exists for almost every $t \in [0,T]$.

Let $t$ be a point where $g'(t)$ exists and $\dot{x}(t) = f(t,x(t))$. At such a point, the time-superdifferential of the composite function collapses to a singleton set (Lemma 1.8 in~\cite{bardi1997optimal}):
\begin{equation}
D^+_t g(t) = \left\{ \frac{d}{dt} J(t,x(t)) \right\}.
\end{equation}
By Proposition 2.7 in~\cite{bardi1997optimal}, we have the inclusion:
\begin{align*}
\left\{ p + q^\top f(t,x(t)) \mid (p,q) \in D^+ J(t,x(t)) \right\} \\ \subseteq \left\{ \frac{d}{dt} J(t,x(t)) \right\}.
\end{align*}
If $D^+ J(t,x(t))$ is empty, the statement holds vacuously. If it is non-empty, let $(p,q)$ be any element in $D^+ J(t,x(t))$. Since the set on the right hand side contains only a single element, the inclusion forces equality:
\begin{equation}
p + q^\top f(t,x(t)) = \frac{d}{dt} J(t,x(t)).
\end{equation}
The result for $D^-$ follows by a similar argument.
\end{proof}

\noindent \textbf{\underline{Discussion: Generalized Chain Rule Implications}}

Lemma~\ref{lem:chain_rule} establishes the generalized chain rule given in Eq.~\eqref{eq: chain rule}. We next discuss the implications of this chain rule, building intuition for how sub/super differentials behave along trajectories.

\paragraph{Trajectories only that move through differentiable regions}

 Since $J$ is Lipschitz, it is differentiable almost everywhere (Rademacher's Theorem). Consequently, if we define$$
S_{J}:=\{(t,x) \in  [0,T] \times  \R^n :  J  \text{ is differentiable} \text{ at } (t,x) \},
$$
then the set $([0,T] \times \R^n /S_J )$ has zero measure. Hence, $\frac{\partial}{\partial t} J(t,x)$ and 
$\nabla_x J(t,x)$ exist for $(t,x) \in [0,T] \times \R^n$ a.e. 

Now, consider a trajectory that stays in differentiable regions of $J$, or only momentarily, with zero measure, cross non-differentiable sets. Specifically, $(t,x(t)) \in S_J$ for all $t \in [0,T]$ a.e. Then $J(t,x(t))$ is differentiable a.e, implying $D^\pm J(t,x(t))$ collapses to a singleton set $\left\{\left(\frac{\partial}{\partial t} J(t,x(t)), \nabla_x J(t,x(t)) \right)  \right\}$ a.e (Lemma 1.8 in~\cite{bardi1997optimal}). For this case Eq.~\eqref{eq: chain rule} reduces to the classical chain rule a.e:
\begin{align} \nonumber 
\frac{d}{dt} J(t,x(t)) = \frac{\partial}{\partial t}J(t,x(t))  & + \nabla_x J(t,x(t)) ^\top f(t,x(t))\\\nonumber 
& \text{ for all } t \in [0,T] \text{ a.e}.
\end{align}

\paragraph{Trajectories that move through non-differentiable regions}

Consider a trajectory that moves along non-differentiable regions of $J$ over a time interval with positive measure, that is $(t,x(t)) \notin S_J$ for $t \in I$ where $I \subset [0,T]$ has positive measure. Along this trajectory the set of sub/superdifferentials, $D^\pm J(t,x(t))$, no longer collapses to a single element. However, the composite function $g(t) := J(t,x(t))$ is absolutely continuous and thus $\frac{d}{dt} J(t,x(t))$ takes a well defined unique singular value a.e. The chain rule given in Eq.~\eqref{eq: chain rule} thus implies that for any $(p_1,q_1),(p_2,q_2) \in D^\pm J(t,x(t))$ we have,
\begin{align*} 
\frac{d}{dt} J(t,x(t)) &= p_1 + q_1^\top f(t,x(t)),\\
\frac{d}{dt} J(t,x(t)) &= p_2 + q_2^\top f(t,x(t)),
\end{align*}
which in turn implies 
\begin{align} \label{eq: orthagonal condition}
\begin{bmatrix} p_1-p_2\\ q_1 -q_2 \end{bmatrix}^\top \begin{bmatrix} 1\\ \dot{x} \end{bmatrix}= p_1-p_2+  (q_1 - q_2)^\top \dot{x} = 0.\end{align} 
That is, the generalized chain rule given Eq.~\eqref{eq: chain rule} implies that the augmented vector field, $\begin{bmatrix} 1 \\ f(t,x) \end{bmatrix}$, is always orthogonal to the {difference} between any two super/sub differentials inside of $D^\pm J(t,x(t))$.

In Example~\ref{ex: non-diff traj} and associated Figure~\ref{fig: moving along the ridge} we give an illustration of how the orthogonality condition in Eq.~\eqref{eq: orthagonal condition} is satisfied. For this example, if $(p,q) \in D^\pm J(t,x(t))$ then $p=0$ and $q^\top \dot{x}=0$, that is $J$ is time invariant and has spatial sub/super derivatives that point in the orthogonal direction of the vector field.

\begin{figure}
\centering
\begin{tikzpicture}[thick, scale=0.55]
\begin{axis}[
view={125}{20}, 
width=12cm, height=10cm,
xmin=-2, xmax=2,
ymin=-2, ymax=2,
zmin=-2, zmax=2.5,
xlabel={$x_1$},
ylabel={$x_2$},
zlabel={$J(x)$},
label style={font=\boldmath}, 
label style={font=\large\bfseries},
tick label style={font=\footnotesize},
axis lines=left,
z buffer=sort,
colormap={ocean}{color=(blue!10) color=(blue!60)},
legend style={
	at={(0.02,0.98)}, 
	anchor=north west, 
	cells={anchor=west},
	font=\small,
	fill=white, 
	fill opacity=0.9,
	draw opacity=0.5 } ]
\addlegendimage{area legend, fill=green!60!black, opacity=0.4}
\addlegendentry{Green regions: Tangent planes}
\addlegendimage{area legend, shade, top color=blue!40, bottom color=blue!40, opacity=0.8}
\addlegendentry{Blue region: $J=|x_1-x_2|$}
\addlegendimage{line legend, red, line width=2pt}
\addlegendentry{Trajectory $x(t)=[t;t]$}
\pgfplotsinvokeforeach{-0.5,-0.3,-0.1,0.1,0.3,0.5}{
	\addplot3[
	surf,
	shader=flat,
	fill=green!60!black,
	opacity=0.20,
	draw=none,
	domain=-2:2,
	domain y=-2:2,
	samples=2,
	forget plot
	]
	{#1*(x-y)};
}
\addplot3[
surf,
shader=interp,
opacity=0.75,
draw=none,
domain=-2:2,
domain y=-2:2,
samples=20,
forget plot
]
{abs(x-y)};
\draw[
red,
line width=2.5pt,
-{Stealth[length=4mm,width=2.5mm]}
]
(axis cs:-2,-2,0.1) -- (axis cs:1.5,1.5,0.1);
\node[
red,
anchor=south west,
font=\bfseries
]
at (axis cs:1.5,0.9,0.9)
{\large $f(x)=[1;1]$};
\draw[
black,
line width=2pt,
-{Stealth[length=4mm,width=2.5mm]}
]
(axis cs:1.5,1.5,0.1) -- (axis cs:1.5,1.5,-1.3);

\node[
black,
anchor=north,
font=\bfseries\small,
align=center
]
at (axis cs:1.5,1.5,-1.55)
{tangential plane\\normals};

\draw[black, thick]
(axis cs:1.5,1.5,-0.15) --
(axis cs:1.62,1.38,-0.15) --
(axis cs:1.62,1.38,0.1);

\draw[
black,
line width=1.5pt,
Stealth-Stealth
]
plot[
domain=-45:45,
samples=5,
variable=\t
] ({1.5 + 0.35*sin(\t)},{1.5 - 0.35*sin(\t)},{-1.3 + 0.35*(1-cos(\t))});

\end{axis}
\end{tikzpicture}
\vspace{-8pt}
\caption{\footnotesize Visualization of Example~\ref{ex: non-diff traj}. The plot shows the surface $\mathcal{S} = \{x \in \R^3: |x_1-x_2|-x_3=0\}$. The black curved arrow represents the set of surface normals $\mathbf{n} = [\alpha,-\alpha,-1]^\top$ with $\alpha \in [-1,1]$, derived from the gradient of the implicit surface equation. The red trajectory moves along the singular ridge (the region of non-differentiability) where the tangent vector is orthogonal to the entire set of subgradients.}
\label{fig: moving along the ridge}
\vspace{-10pt}
\end{figure}

\begin{ex} \label{ex: non-diff traj}
Consider the function $J(x)=|x_1-x_2|$ and the trajectory $x(t)=[t, t]^\top$ with dynamics $\dot{x}=[1, 1]^\top$. Along this trajectory, $J(x(t))=0$ and thus $\frac{d}{dt} J(x(t))=0$.
By Theorem 23.9 in~\cite{Rockafellar_1970} and Proposition 4.7 in~\cite{bardi1997optimal}, the sub differential set is $D^\pm J(x(t)) = \{ [\alpha,-\alpha]^\top: \alpha \in [-1,1] \}$.
We can then verify the chain rule of Lemma~\ref{lem:chain_rule}:
\[
q^\top \dot{x}(t)  = [\alpha,-\alpha]^\top [1, 1] = \alpha - \alpha = 0 = \frac{d}{dt} J(x(t))
\]
for any sub differential $q \in D^\pm J(x(t))$ and $t>0$. This geometry is illustrated in Figure~\ref{fig: moving along the ridge}.
\end{ex}

\paragraph{Lipschitz functions with empty sub/super differential sets}
As we will see in Example~\ref{ex: saddle point empty}, there exist Lipschitz continuous functions that have no sub/super differentials at certain points. That is it is possible that ${D^\pm J(t,x(t))} = \emptyset$. For these functions, we are no longer able to apply the chain rule given in Eq.~\eqref{eq: chain rule} since this equation only holds whenever there exists $(p,q) \in D^\pm 	J (t,x(t))$ implying $D^\pm 	J (t,x(t)) \ne \emptyset$.

\begin{ex} \label{ex: saddle point empty} Consider the saddle function $J(t,x)=|x_1-t|-|x_2-t|$ and the trajectory $x(t)=[t, t]^\top$ with dynamics $\dot{x}=[1, 1]^\top$. Along this trajectory, $J(t,x(t))=0 $ and thus  $\frac{d}{dt}J(t,x(t))=0$. 

However, at the point $x(t)$, $J$ looks like a ``V'' shape along the $x_1$-axis and an inverted ``V'' along the $x_2$-axis, creating an ``X'' shape cross-section. Unlike in Fig.~\ref{fig: moving along the ridge}, it is impossible to place a tangent plane above or below the centre of this "X" shape without it cutting through the surface of $J$. Any plane attempting to touch the point from below will inevitably intersect the inverted ``V'' arms, and any plane touching from above will intersect the upright ``V''. Since no tangent plane can sit perfectly above or below the surface it follows $D^\pm 	J (t,x(t)) = \emptyset$.

Consequently, while Lemma~\ref{lem:chain_rule} holds vacuously (whenever there exists $(p,q) \in D^\pm 	J (t,x(t))$ implying $D^\pm 	J (t,x(t)) \ne \emptyset$), the explicit chain rule given in Eq.~\eqref{eq: chain rule} cannot be instantiated along this trajectory for any $t \ge 0$. There exists no $(p,q) \in D^\pm J(t,x(t))$ such that $ p+ q^\top [1,1]^\top =0= \frac{d}{dt} J(t,x(t))$. 
\end{ex}

Note, although saddle functions, such as the function given in Example~\ref{ex: saddle point empty}, are Lipschitz continuous, Definition~\ref{def: admis VF} implicitly restricts these type of functions from being used as candidate value functions. Definition~\ref{def: admis VF}, ensures that all admissible candidate value functions, $J$, satisfy the property that ${D^\pm J(t,x)} \ne \emptyset$ for all $(t,x) \in [0,T] \times  \R^n $, and hence we can always apply the chain rule in Eq.~\eqref{eq: chain rule} to admissible candidate value functions. This is because the definition requires ${D^\pm J(t,x)} = \partial_C J(t,x)$  and $\partial_C J(t,x) \ne \emptyset$ (Remark 5.1~\cite{zhou1993verification} or Page 196~\cite{clarke2013functional}).

\begin{figure}
\centering
\begin{tikzpicture}[thick, scale=0.55]
\begin{axis}[
view={150}{25}, 
width=12cm, height=10cm,
xmin=-2, xmax=2,
ymin=-2, ymax=2,
zmin=-2.5, zmax=2.5,
xlabel={$x_1$},
ylabel={$x_2$},
zlabel={$J(0,x)$},
label style={font=\large\bfseries},
tick label style={font=\footnotesize},
axis lines=left,
z buffer=sort,
colormap={saddle}{color=(red!80!black) color=(blue!80!black) color=(red!80!black)},
legend style={
	at={(0.02,0.98)},
	anchor=north west,
	cells={anchor=west},
	font=\small,
	fill=white,
	fill opacity=0.9,
	draw opacity=0.4 } ]

\addlegendimage{area legend, shade, top color=red!70!black, bottom color=red!70!black, middle color=blue!70!black, opacity=0.4}
\addlegendentry{Saddle Surface: $J=|x_1-t|-|x_2-t|$}
\addlegendimage{only marks, mark=*,
	mark options={scale=1.5, fill=red},
	color=red}
\addlegendentry{Trajectory $x(t)=[t;t]$}
\addlegendimage{only marks, mark=*, mark options={fill=red}, color=red}

\addplot3[
surf,
shader=interp,
opacity=0.4,
draw=none,
domain=-2:2,
domain y=-2:2,
samples=20, 
forget plot
]
{abs(x) - abs(y)};



\addplot3[
only marks,
mark=*,
mark options={scale=1.5, fill=red},
color=red,
forget plot
] coordinates {(0,0,0)};

\node[
red,
anchor=south,
font=\bfseries\small,
align=center,
yshift=5pt
] at (axis cs:0,-0.5,	1) { \Large $D^\pm J(t,[t,t]^\top)=\emptyset$};

\end{axis}
\end{tikzpicture}
\caption{\footnotesize Visualization of Example~\ref{ex: saddle point empty} involving the saddle function $J(t,x) = |x_1-t| - |x_2-t|$ at $t=0$. The red trajectory $x(t)=[t,t]^\top$ moves through the non-differentiable region where both sub/super differentials do not exist.	}
\label{fig: saddle point failure}
\end{figure}
\section*{Acknowledgements}
This research was supported by the NSF Grant 2606245.
\else
\section{appendix}
{\begin{lemma} \label{lem: ess sup bound}
		Let $E : [0,T] \times \mathbb{R}^n \to \mathbb{R}$ be a lower semicontinuous function. For any measurable function $x: [0,T] \to B_R(0) \subset \mathbb{R}^n$, for some $R>0$, the following inequality holds:
		\begin{equation}
			\operatorname*{ess\,sup}_{t \in [0,T]} E(t, x(t)) \le \operatorname*{ess\,sup}_{(t,y) \in [0,T] \times B_R(0)} E(t, y).
		\end{equation}
	\end{lemma}
    \begin{proof}
        Since $E$ is lower semicontinuous, its strict superlevel sets are open, meaning any exceedance of the essential supremum along $(t,x(t))$ would contradictorily require an exceedance on a non-empty open euclidean set which by definition has strictly positive Lebesgue measure (see Arxiv version for more details~\cite{jones2025bounding}).
    \end{proof}

\begin{lemma} \label{lem:chain_rule}
Consider functions $J \in \text{Lip}([0,T] \times \mathbb{R}^n , \mathbb{R})$ and $x \in {\text{AC}}([0,T], \mathbb{R}^n)$ satisfying $\dot{x}(t) = f(t,x(t))$ almost everywhere on $[0,T]$.
Then, for almost every $t \in [0,T]$ and any $(p,q) \in D^\pm J(t,x(t))$, we have
\begin{align} \label{eq: chain rule}
\frac{d}{dt} J(t,x(t)) &= p + q^\top f(t,x(t)).
\end{align}
\end{lemma}
\begin{proof}
	The composite function $t \mapsto J(t,x(t))$ is absolutely continuous, guaranteeing its classical derivative exists almost everywhere. At points of differentiability, the time-superdifferential reduces to the singleton $\{ \frac{d}{dt} J(t,x(t)) \}$ \cite[Lemma 1.8]{bardi1997optimal}, which, combined with the chain rule inclusion for superdifferentials \cite[Proposition 2.7]{bardi1997optimal}, directly implies the desired equality (see Arxiv version for more details~\cite{jones2025bounding}).
\end{proof} }

\section*{Acknowledgements}
This research was supported by the NSF Grant 2606245.

\fi

\end{document}